\documentclass[11pt]{amsart}

\usepackage[top=24mm, bottom=24mm, left=24mm, right=24mm]{geometry}
\usepackage{color}
\usepackage{graphicx}
\usepackage{bm}

\newcommand{\ngauss}{q}
\newcommand{\ncheb}{p}

\newcommand{\pvct}[1]{\bm{#1}}
\newcommand{\vct}[1]{\bm{\mathsf{#1}}}

\newcommand{\pxx}{\pvct{x}}
\newcommand{\pyy}{\pvct{y}}
\newcommand{\pzz}{\pvct{z}}
\newcommand{\uu}{\vct{u}}
\newcommand{\vv}{\vct{v}}
\newcommand{\ww}{\vct{w}}

\renewcommand{\AA}{\mtx{A}}

\newcommand{\TT}{\mtx{T}}
\newcommand{\UU}{\mtx{U}}
\newcommand{\VV}{\mtx{V}}
\newcommand{\WW}{\mtx{W}}

\newcommand{\mtx}[1]{\bm{\mathsf{#1}}}

\newtheorem{remark}{Remark}
\theoremstyle{definition}

\newcommand{\pgnotate}[1]{}

\newcommand{\lsp}{\vspace{3mm}}

\begin{document}

\begin{center}
\textbf{\large{The Hierarchical Poincar\'e-Steklov (HPS) solver for elliptic PDEs: A tutorial}}
\label{chapter5:spectralcomposite}

\lsp

\textit{P.G.~Martinsson, Dept.~of Applied Mathematics, University of Colorado at Boulder}

\lsp

\textit{June 3, 2015}

\lsp

\begin{minipage}{130mm}\small
\textbf{Abstract:}
A numerical method for variable
coefficient elliptic problems on two dimensional domains is
described. The method is based on high-order spectral approximations
and is designed for problems with smooth solutions. The resulting
system of linear equations is solved using a direct solver with
$O(N^{1.5})$ complexity for the pre-computation and $O(N \log N)$
complexity for the solve.  The fact that the solver is direct is a
principal feature of the scheme, and makes it particularly well suited
to solving problems for which iterative solvers struggle; in particular
for problems with highly oscillatory solutions. This note is intended
as a tutorial description of the scheme, and draws heavily on previously
published material.
\end{minipage}
\end{center}

\section{Introduction}

\subsection{Problem formulation and outline of solution strategy}
This note describes a direct solver for elliptic PDEs with
variable coefficients, such as, e.g.,
\begin{equation}
\label{eq5:basic}
\left\{\begin{aligned}
\mbox{}[Au](\pxx) =&\ g(\pxx),\qquad &\pxx \in \Omega,\\
   u(\pxx) =&\ f(\pxx),\qquad &\pxx \in \Gamma,
\end{aligned}\right.
\end{equation}
where $A$ is a variable coefficient elliptic differential operator
\begin{multline}
\label{eq5:defA}
[Au](\pxx) = -c_{11}(\pxx)[\partial_{1}^{2}u](\pxx)
-2c_{12}(\pxx)[\partial_{1}\partial_{2}u](\pxx)
-c_{22}(\pxx)[\partial_{2}^{2}u](\pxx)\\
+c_{1}(\pxx)[\partial_{1}u](\pxx)
+c_{2}(\pxx)[\partial_{2}u](\pxx)
+c(\pxx)\,u(\pxx),
\end{multline}
where $\Omega$ is a box in the plane with boundary $\Gamma = \partial \Omega$,
where all coefficient functions ($c$, $c_{i}$, $c_{ij}$) are smooth,
and where $f$ and $g$ are given functions.
(For generalizations, see Section \ref{sec5:generalization}.)
The solver is structured as follows:
\begin{enumerate}
\item The domain is first tessellated into a hierarchical tree
of patches. For each patch on the finest level, a reduced model
that we call a ``proxy''
that represents its internal structure is computed. The proxy takes
the form of a dense matrix (that may be stored in a data sparse format),
and is for the small patches computed by brute force.
\item The larger patches are processed
in an upwards pass through the tree. For each patch, its proxy matrix is
formed by merging the proxies of its children.
\item Once the proxies for all patches have been computed, a solution
to the PDE can be computed via a downwards pass through the tree.
This step is typically very fast.
\end{enumerate}
We observe that this pattern is similar to the classical
nested dissection method of George \cite{george_1973}, with
a large subsequent literature on ``multifrontal solver,'' see,
e.g., \cite{1989_directbook_duff,2006_davis_directsolverbook}
and the references therein.

The techniques described in this note are drawn from
\cite{2012_martinsson_composite_orig,2012_spectralcomposite,2013_martinsson_ItI,2013_martinsson_DtN_linearcomplexity},
which in turn is inspired by earlier work on multidomain spectral methods, see,
e.g., \cite{2003_pfeiffer_spectralmultidomain,1999_hesthaven_spectralcollocation} and the references therein.

\subsection{Dirichlet-to-Neumann maps}
In this section, the internal structure of a patch is represented
by computing its \textit{Dirichlet-to-Neumann}, or ``DtN,'' map. To explain
what this map does, first observe that for a given boundary function
$f$, the BVP (\ref{eq5:basic}) typically has a unique solution $\phi$
(unless the operator $A$ happens to have a non-trivial null-space, see
Remark \ref{remark:resonance}). Now simply form the boundary function $h$
that gives the normal derivative of the solution,
$$
h(\pxx) = \phi_{n}(\pxx),\qquad \pxx \in \Gamma,
$$
where $\phi_{n}$ is the outwards pointing normal derivative.
The process for constructing the function $h$ from $f$ is linear, and
we write it as
$$
h = T\,f.
$$
Or, equivalently,
$$
T\,\colon\,\phi|_{\Gamma} \mapsto \phi_{n}|_{\Gamma},
\qquad\qquad \mbox{where}\ \phi\ \mbox{satisfies}\
A\phi = 0.
$$

In general, the map $T$ is a slightly unpleasant object; it behaves
as a differentiation operator, and it has complicated singular behavior
near the corners of $\Gamma$. A key observation is that in the
present context, all these difficulties can be
ignored since we limit attention to functions that are smooth. In a
sense, we only need to accurately represent the projection of the
``true'' operator $T$ onto a space of smooth functions (that in particular
do not have any corner singularities).

We represent boundary functions by tabulating their values at
interpolation nodes on the edges of the boxes. For instance,
for a leaf box, we place $q$ Gaussian nodes on each side
(for say $q=20$), which means that the functions $f$ and $h$
are represented by vectors $\vct{f},\vct{h} \in \mathbb{R}^{4q}$
and the discrete approximation to $T$ is a $4q \times 4q$
matrix $\mtx{T}$. The technique for computing $\mtx{T}$ for a
leaf box is described in Section \ref{sec5:leaf}. For a parent
box $\tau$ with children $\alpha$ and $\beta$, there
is a technique for computing $\mtx{T}^{\tau}$ from the
matrices $\mtx{T}^{\alpha}$ and $\mtx{T}^{\beta}$
that is described in Section \ref{sec5:merge}. In essence, the
idea is simply to enforce continuity of both potentials and
fluxes across the edge that is shared by $\Omega_{\alpha}$
and $\Omega_{\beta}$.

\begin{remark}
\label{remark:resonance}
For a general BVP like (\ref{eq5:basic}), the DtN operator need not
exist. For example, suppose $-k^{2}$ is an eigenvalue of
$\Delta$ with zero Dirichlet data. Then the operator $Au = \Delta u + k^{2}u$
clearly has a non-trivial null-space. However, this situation is
in some sense ``unusual'' (since the spectrum of an elliptic PDO
on a bounded domain typically is discrete), and it turns out that
one can for the most part completely ignore this complication and
assume that the DtN operator always exists. Note for instance that
if $A$ is coercive (e.g.~if $A = -\Delta$), then the DtN is guaranteed
to exist for any bounded domain. For problems for which resonances
do present numerical problems, there is a variation of the proposed
method that is rock-solid stable. The idea is to build a hierarchy
of so called \textit{impedance maps} (instead of DtN maps). These are
cousins of the DtN that are defined by
$$
R\,\colon\,
(\phi + \mathrm{i}\phi_{n})|_{\Gamma} \mapsto
(\phi - \mathrm{i}\phi_{n})|_{\Gamma},
\qquad\qquad \mbox{where}\ \phi\ \mbox{satisfies}\
A\phi = 0.
$$
The map $R$ always exists, and is moreover a unitary operator. This
construction was proposed by Alex Barnett of Dartmouth College
\cite{2013_martinsson_ItI}.
\end{remark}

\subsection{Complexity of the direct solver}

The asymptotic complexity of the solver described in this note
is exactly the same as that for classical nested dissection \cite{george_1973}.
For a domain with $N$ interior discretization nodes, the pre-computation
(the upwards pass) costs $O(N^{1.5})$ operations, and then the
solve (the downwards pass) costs $O(N\log N)$ operations,
see \cite[Sec.~5.2]{2012_spectralcomposite}.

Optimal $O(N)$ complexity can be achieved for both the pre-computation
and the solve stage when the Green's function of the BVP is non-oscillatory.
In this case, the matrices $\mtx{T}^{\tau}$
that approximate the DtN operators have enough internal structure that they
can be well represented using a data-sparse format such as, e.g., the
$\mathcal{H}$-matrix framework of Hackbusch and co-workers
\cite{hackbusch,2002_hackbusch_H2,2008_bebendorf_book,2010_borm_book},
see \cite{2013_martinsson_DtN_linearcomplexity}.

\subsection{Generalizations}
\label{sec5:generalization}
For notational simpliticy, this note treats only the simple Dirichlet
problem (\ref{eq5:basic}). The scheme can with trivial modifications be
applied to more general elliptic operators
coupled with Dirichlet, Neumann, or mixed boundary data. It has for
instance been successfully tested on convection-diffusion problems
that are strongly dominated (by a factor of $10^{4}$) by the
convection term. It has also been tested on vibration problems
modeled by a variable coefficient Helmholtz equation, and has
proven capable of solving problems on domains of size $200 \times 200$
wavelengths or more on an office laptop computer. The solutions
are computed to seven correct digits. Extension to more general
domains is done via parameter maps to a rectangle, or union of
rectangles. See \cite{2012_spectralcomposite} for details.

\subsection{Outline}

To keep the presentation uncluttered, we start by describing a direct
solver for (\ref{eq5:basic}) for the case of no body-load, $g=0$, and
with $\Omega = [0,1]^{2}$ the unit square.
Section \ref{sec5:discretization} introduces the discretization, and
how we represent the approximation to the solution $u$ for (\ref{eq5:basic}).
Section \ref{sec5:leaf} describes the how to compute the DtN operator for a leaf in the tree.
Section \ref{sec5:merge} describes the merge process for how to take the DtN operators for
two touching boxes, and computing the DtN operator for their union.
Section \ref{sec5:hierarchy} describes the full hierarchical scheme.
Section \ref{sec5:body} describes how to solve a problem with body
loads.

\section{Discretization}
\label{sec5:discretization}

Partition the domain $\Omega$ into a collection of square (or
possibly rectangular) boxes, called \textit{leaf boxes}.
On the edges of each leaf, place $\ngauss$ Gaussian interpolation
points. The size of the leaf boxes, and the parameter $\ngauss$
should be chosen so that any potential solution $u$ of
(\ref{eq5:basic}), as well as its first and second derivatives,
can be accurately interpolated from their values at these points
($\ngauss=20$ is often a good choice).
Let $\{\pxx_{k}\}_{k=1}^{N}$ denote the collection of interpolation
points on all boundaries.

\begin{figure}
\includegraphics[width=\textwidth]{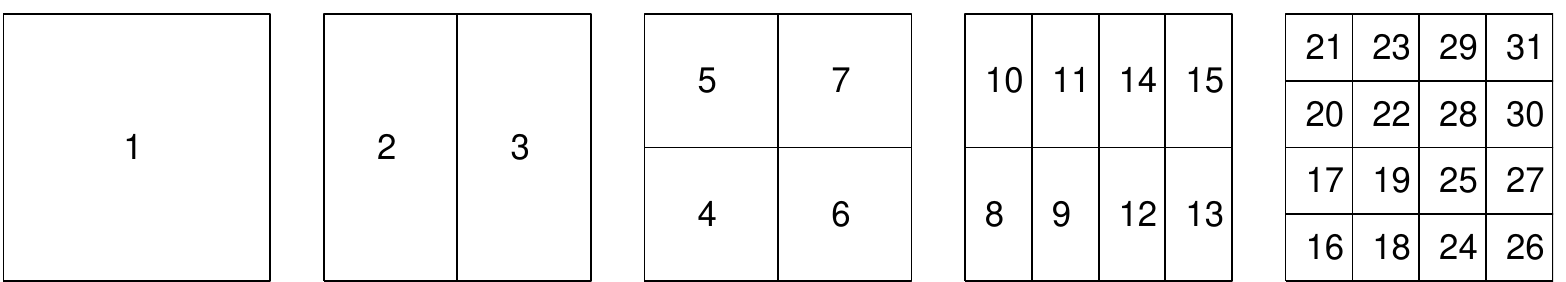}
\caption{The square domain $\Omega$ is split into $4 \times 4$ leaf boxes.
These are then gathered into a binary tree of successively larger boxes
as described in Section \ref{sec5:thealgorithm}. One possible enumeration
of the boxes in the tree is shown, but note that the only restriction is
that if box $\tau$ is the parent of box $\sigma$, then $\tau < \sigma$.}
\label{fig5:tree_numbering}
\end{figure}

Next construct a binary tree on the collection of leaf boxes by
hierarchically merging them, making sure that all boxes on
the same level are roughly of the same size, cf.~Figure
\ref{fig5:tree_numbering}.  The boxes should be ordered so
that if $\tau$ is a parent of a box $\sigma$, then $\tau < \sigma$. We
also assume that the root of the tree (i.e.~the full box $\Omega$) has
index $\tau=1$. We let $\Omega_{\tau}$ denote the domain associated with box $\tau$.

With each box $\tau$, we define two index vectors $I_{\rm i}^{\tau}$ and $I_{\rm e}^{\tau}$ as follows:

\vspace{1mm}

\begin{itemize}
\item[$I_{\rm e}^{\tau}$] A list of all \textit{exterior} nodes of $\tau$.
In other words,
$k \in I_{\rm e}^{\tau}$ iff $\pxx_{k}$ lies on the boundary of $\Omega_{\tau}$.

\vspace{1mm}

\item[$I_{\rm i}^{\tau}$]
For a parent $\tau$, $I_{\rm i}^{\tau}$ is a list of all its \textit{interior} nodes
that are not interior nodes of its children.
For a leaf $\tau$, $I_{\rm i}^{\tau}$ is empty.
\end{itemize}

\vspace{1mm}

\noindent
Let $\uu \in \mathbb{R}^{N}$ denote a vector holding approximations
to the values of $u$ of (\ref{eq5:basic}), in other words,
$$
\uu(k) \approx u(\pxx_{k}).
$$
Finally, let $\vv \in \mathbb{R}^{N}$ denote a vector holding approximations to
the boundary fluxes of the solution $u$ of (\ref{eq5:basic}), in other words
$$
\vv(k) \approx \left\{
\begin{aligned}
\partial_{2}u(\pxx_{k}),\qquad&\mbox{when }\pxx_{j}\mbox{ lies on a horizontal edge,}\\
\partial_{1}u(\pxx_{k}),\qquad&\mbox{when }\pxx_{j}\mbox{ lies on a vertical edge.}
\end{aligned}\right.
$$
Note the sign convention for the normal derivatives: we use a global frame
of reference, as opposed to distinguishing between outwards and inwards pointing
normal derivatives. This is a deliberate choice to avoid problems with signs when
matching fluxes of touching boxes.


\section{Constructing the Dirichlet-to-Neumann map for a leaf}
\label{sec5:leaf}

This section describes a spectral method for computing a discrete approximation to the DtN
map $T^{\tau}$ associated with a leaf box $\Omega_{\tau}$. In other words, if $u$ is a
solution of (\ref{eq5:basic}), we seek a matrix $\mtx{T}^{\tau}$ of size $4\ngauss \times 4\ngauss$
such that
\begin{equation}
\label{eq5:Ttau}
\begin{array}{ccccc}
\vv(I_{\rm e}^{\tau}) &\approx& \mtx{T}^{\tau} & \uu(I_{\rm e}^{\tau}).\\
\mbox{\textit{Neumann data}} && \mbox{\textit{DtN map}} & \mbox{\textit{Dirichlet data}}
\end{array}
\end{equation}
Conceptually, we proceed as follows: Given a vector $\uu(I_{\rm e}^{\tau}) \in \mathbb{R}^{4q}$ specifying the
solution $u$ on the boundary of $\Omega_{\tau}$, form for each side the unique polynomial
of degree at  most $\ngauss-1$ that interpolates the $\ngauss$ specified values of $u$. This yields Dirichlet
boundary data on $\Omega_{\tau}$ in the form of four polynomials. Solve the restriction
of (\ref{eq5:basic}) to $\Omega_{\tau}$ for the specified boundary data using a spectral method
on a local tensor product grid of $\ncheb \times \ncheb$ \textit{Chebyshev nodes} (typically,
we choose $\ncheb=\ngauss+1$). The vector
$\vv(I_{\rm e}^{\tau})$ is obtained by spectral differentiation of the local solution, and
then re-tabulating the boundary fluxes to the Gaussian nodes in $\{\pxx_{k}\}_{k \in I_{\rm e}^{\tau}}$.

We give details of the construction in Section \ref{sec5:const}, but as a preliminary step,
we first review a classical spectral collocation method for the local solve in Section \ref{sec5:prelim}

\begin{remark}
\label{re:chebvsgauss}
Chebyshev nodes are ideal for the leaf computations, and it is in principle also possible
to use Chebyshev nodes to represent all boundary-to-boundary ``solution operators'' such
as, e.g., $\mtx{T}^{\tau}$ (indeed, this was the approach taken in the
first implementation of the proposed method \cite{2012_spectralcomposite}).
However, there are at least two substantial benefits to using Gaussian nodes that justify
the trouble to retabulate the operators. First, the procedure for merging boundary operators
defined for neighboring boxes is much cleaner and involves less bookkeeping since the
Gaussian nodes do not include the corner nodes. (Contrast Section 4 of \cite{2012_spectralcomposite}
with Section \ref{sec5:merge}.)  Second, and more importantly, the use of the Gaussian
nodes allows for interpolation between different discretizations. Thus the method can
easily be extended to have local refinement when necessary, see Remark \ref{remark:refine}.
%
\end{remark}

\subsection{Spectral discretization}
\label{sec5:prelim}
Let $\Omega_{\tau}$ denote a rectangular subset of $\Omega$ with boundary $\Gamma_{\tau}$, and consider
the local Dirichlet problem with the boundary data given by a ``dummy'' function $\psi$
\begin{equation}
\label{eq5:differential}
\left\{\begin{aligned}
\mbox{}[Au](\pxx) =&\ 0,\qquad &\pxx \in \Omega_{\tau},\\
         u(\pxx) =&\ \psi(\pxx),\qquad &\pxx \in \Gamma_{\tau},
\end{aligned}\right.
\end{equation}
where the elliptic operator $A$ is defined by (\ref{eq5:defA}).
We will construct an approximate solution to (\ref{eq5:differential})
using a classical spectral collocation method described in, e.g.,
Trefethen \cite{2000_trefethen_spectral_matlab}:
First, pick a small integer $\ncheb$ and let
$\{\pzz_{k}\}_{k=1}^{\ncheb^{2}}$
denote the nodes in a tensor product grid of $\ncheb\times \ncheb$ Chebyshev
nodes on $\Omega_{\tau}$.
Let $\mtx{D}^{(1)}$ and $\mtx{D}^{(2)}$ denote spectral differentiation matrices corresponding to
the operators $\partial/\partial x_1$ and $\partial/\partial x_2$, respectively.
The operator (\ref{eq5:defA}) is then locally approximated via the $\ncheb^{2} \times \ncheb^{2}$ matrix
\begin{equation}
\AA =
-\mtx{C}_{11}\bigl(\mtx{D}^{(1)}\bigr)^{2}
-2\mtx{C}_{12}\mtx{D}^{(1)}\mtx{D}^{(2)}
-\mtx{C}_{22}\bigl(\mtx{D}^{(2)}\bigr)^{2}
+\mtx{C}_{1}\mtx{D}^{(1)}
+\mtx{C}_{2}\mtx{D}^{(2)}
+\mtx{C},
\end{equation}
where $\mtx{C}_{11}$ is the diagonal matrix with diagonal entries $\{c_{11}(\pzz_{k})\}_{k=1}^{\ncheb^{2}}$,
and the other matrices $\mtx{C}_{ij}$, $\mtx{C}_{i}$, $\mtx{C}$ are defined analogously.

\begin{figure}
\begin{tabular}{ccc}
\includegraphics[height=50mm]{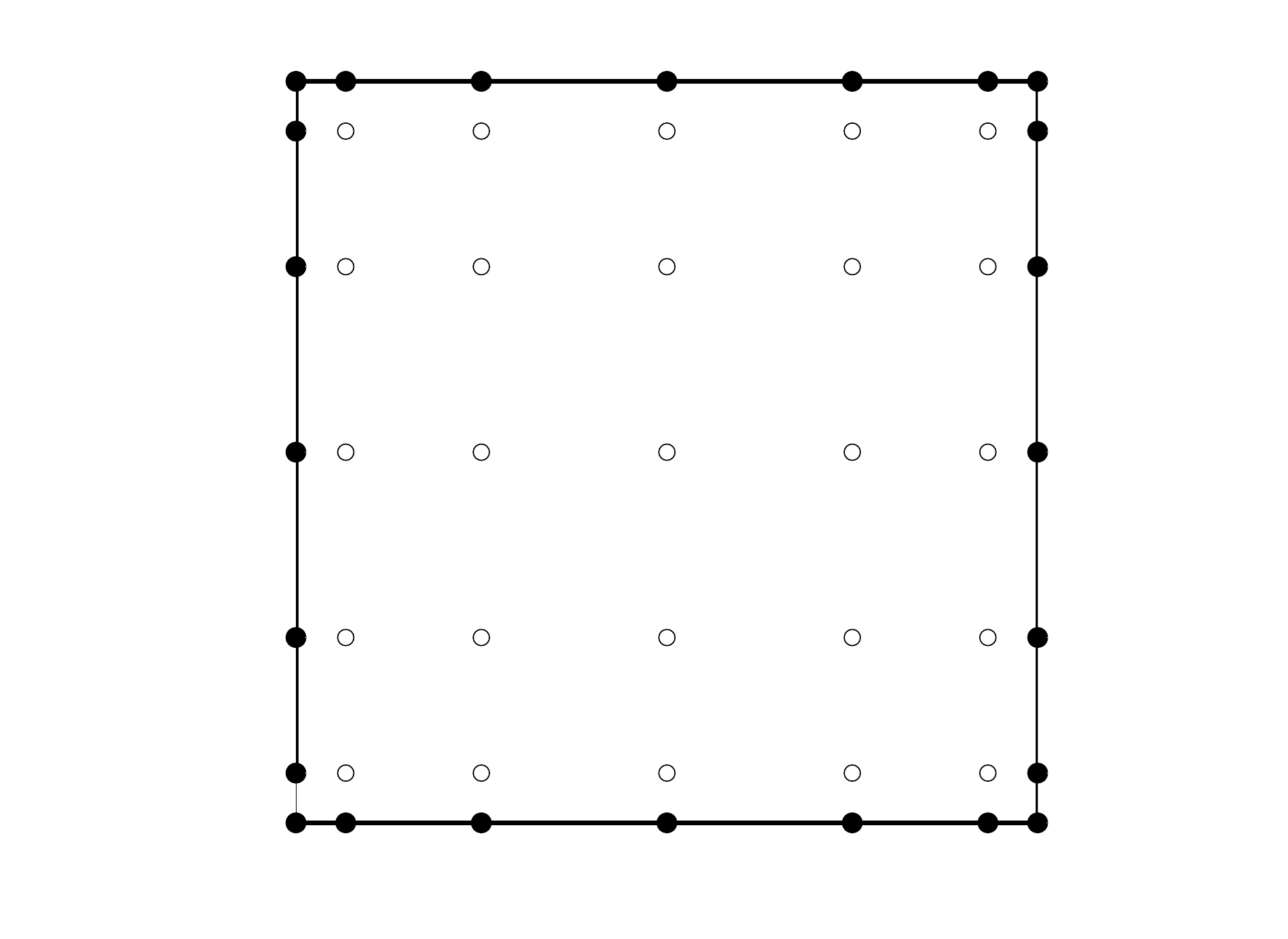}
& \mbox{}\hspace{10mm}\mbox{} &
\includegraphics[height=50mm]{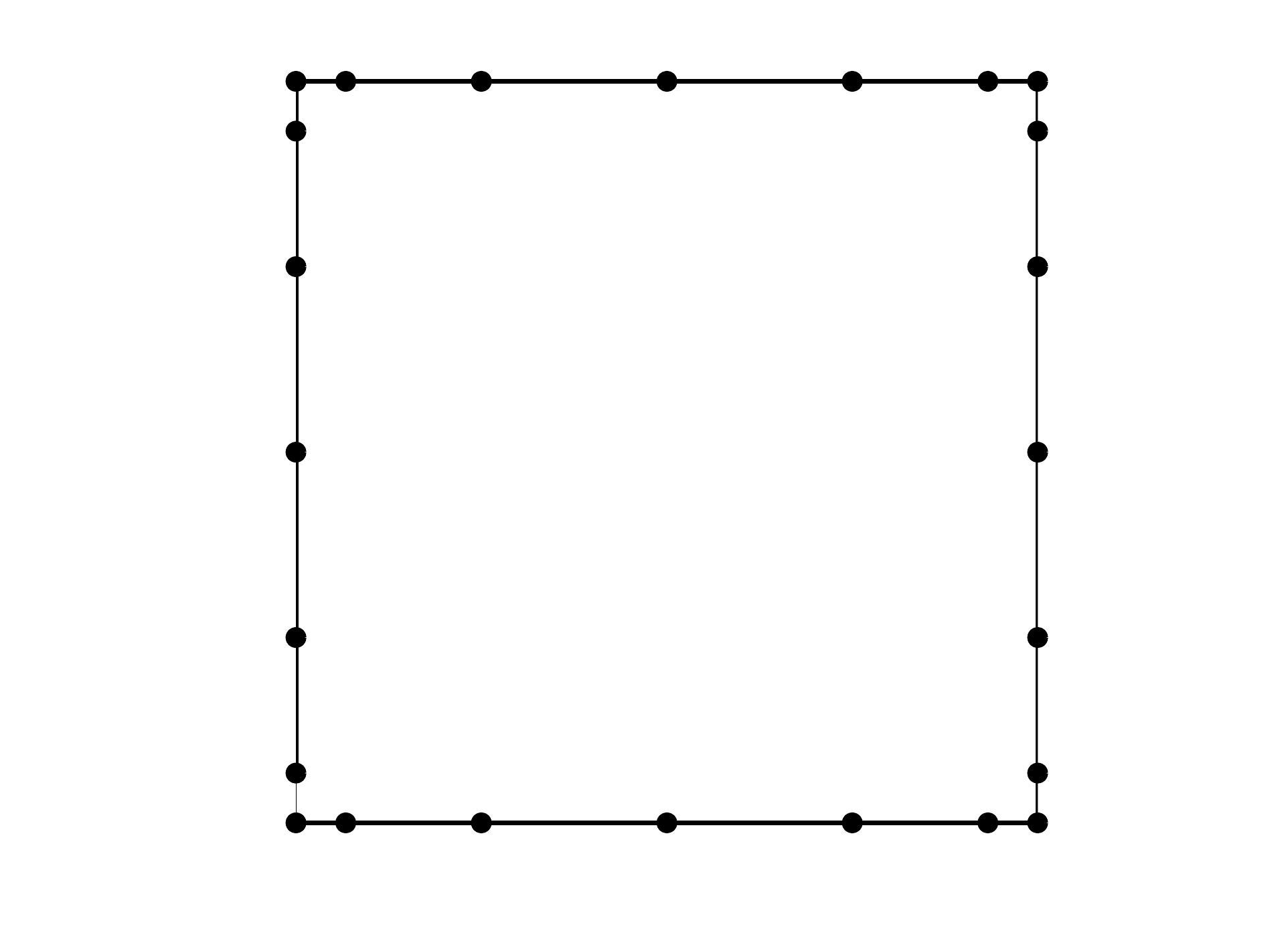}
\\
(a) && (b)
\end{tabular}
\caption{Notation for the leaf computation in Section \ref{sec5:leaf}.
(a) A leaf before elimination of interior (white) nodes.
(b) A leaf after elimination of interior nodes.}
\label{fig5:leaf}
\end{figure}

Let $\ww \in \mathbb{R}^{\ncheb^{2}}$ denote a vector holding the desired
approximate solution of (\ref{eq5:differential}). We populate all entries corresponding
to boundary nodes with the Dirichlet data from $\psi$, and then enforce a spectral
collocation condition at the interior nodes. To formalize, let us partition the index set
$$
\{1,\,2,\,\dots,\,\ncheb^{2}\} = J_{\rm e}\cup J_{\rm i}
$$
in such a way that $J_{\rm e}$ contains the $4(\ncheb-1)$ nodes on the boundary of $\Omega_{\tau}$, and
$J_{\rm i}$ denotes the set of $(\ncheb-2)^{2}$ interior nodes, see Figure \ref{fig5:leaf}(a).
Then partition the vector $\ww$ into two parts corresponding to internal and
exterior nodes via
$$
\ww_{\rm i} = \ww(J_{\rm i}),\quad
\ww_{\rm e} = \ww(J_{\rm e}).
$$
Analogously, partition $\AA$ into four parts via
$$
\mtx{A}_{\rm i,i} = \mtx{A}(J_{\rm i},J_{\rm i}),\quad
\mtx{A}_{\rm i,e} = \mtx{A}(J_{\rm i},J_{\rm e}),\quad
\mtx{A}_{\rm e,i} = \mtx{A}(J_{\rm e},J_{\rm i}),\quad
\mtx{A}_{\rm e,e} = \mtx{A}(J_{\rm e},J_{\rm e}).
$$
The potential at the exterior nodes is now given directly from the
boundary condition:
$$
\ww_{\rm e} = \left[\psi(\pzz_{k})\right]_{k \in J_{\rm e}}.
$$
For the internal nodes, we enforce the PDE (\ref{eq5:differential}) via direct collocation:
\begin{equation}
\label{eq5:A_interior}
\AA_{\rm i,i}\,\ww_{\rm i} +
\AA_{\rm i,e}\,\ww_{\rm e} = \vct{0}.
\end{equation}
Solving (\ref{eq5:A_interior}) for $\ww_{\rm i}$, we find
\begin{equation}
\label{eq5:leaf_solve}
\ww_{\rm i} = -\AA_{\rm i,i}^{-1}\,\AA_{\rm i,e}\,\ww_{\rm e},
\end{equation}

\subsection{Constructing the approximate DtN}
\label{sec5:const}
Now that we know how to approximately solve the local Dirichlet problem
(\ref{eq5:differential}) via a local spectral method, we can
build a matrix $\mtx{T}^{\tau}$ such that (\ref{eq5:Ttau})
holds to high accuracy. The starting point is a vector $\uu(I_{\tau}) \in \mathbb{R}^{4\ngauss}$ of
tabulated potential values on the boundary of $\Omega_{\tau}$. We will
construct the vector $\vv(I_{\tau}) \in \mathbb{R}^{4\ngauss}$ via four linear maps. The combination
of these maps is the matrix $\mtx{T}^{\tau}$.

\lsp

\noindent
\textbf{\textit{Step 1 --- re-tabulation from Gaussian nodes to Chebyshev nodes:}}
For each side of $\Omega_{\tau}$, form the unique interpolating polynomial
of degree at most $\ngauss-1$ that interpolates the $\ngauss$ potential values on
that side specified by $\uu(I_{\rm e}^{\tau})$. Now evaluate these polynomials at the
boundary nodes of a $\ncheb \times \ncheb$ Chebyshev grid on $\Omega_{\tau}$.
Observe that for a corner node, we may in the general case get conflicts. For
instance, the potential at the south-west corner may get one value from extrapolation
of potential values on the south border, and one value from extrapolation of the
potential values on the west border. We resolve such conflicts by assigning the corner
node the average of the two possibly different values. (In practice, essentially no
error occurs since we know that the vector $\uu(I_{\rm e}^{\tau})$ tabulates an
underlying function that is continuous at the corner.)

\lsp

\noindent
\textbf{\textit{Step 2 --- spectral solve:}}
Step 1 populates the boundary nodes of the $\ncheb \times \ncheb$ Chebyshev grid with
Dirichlet data. Now determine the potential at all interior points on the Chebyshev
grid by executing a local spectral solve, cf.~equation (\ref{eq5:leaf_solve}).

\lsp

\noindent
\textbf{\textit{Step 3 --- spectral differentiation:}}
After Step 2, the potential is known at all nodes on the local Chebyshev grid.
Now perform spectral differentiation to evaluate approximations to $\partial u/\partial x_{2}$
for the Chebyshev nodes on the two horizontal sides, and $\partial u/\partial x_{1}$ for the
Chebyshev nodes on the two vertical sides.

\lsp

\noindent
\textbf{\textit{Step 4 --- re-tabulation from the Chebyshev nodes back to Gaussian nodes:}}
After Step 3, the boundary fluxes on $\partial \Omega_{\tau}$ are specified by four polynomials
of degree $\ncheb-1$ (specified via tabulation on the relevant Chebyshev nodes). Now simply
evaluate these polynomials at the Gaussian nodes on each side to obtain the vector $\vv(I_{\rm e}^{\tau})$.

\lsp

Putting everything together, we find that the matrix $\mtx{T}^{\tau}$ is given as a product of
four matrices
$$
\begin{array}{ccccccccc}
\mtx{T}^{\tau} &=& \mtx{L}_{4} &\circ& \mtx{L}_{3} &\circ& \mtx{L}_{2} &\circ& \mtx{L}_{1} \\
4\ngauss \times 4\ngauss
&&
4\ngauss \times 4\ncheb
&&
4\ncheb \times \ncheb^{2}
&&
\ncheb^{2} \times 4(\ncheb-1)
&&
4(\ncheb-1) \times 4\ngauss
\end{array}
$$
where $\mtx{L}_{i}$ is the linear transform corresponding to ``Step $i$'' above. Observe that
many of these transforms are far from dense, for instance, $\mtx{L}_{1}$ and $\mtx{L}_{4}$
are $4\times 4$ block matrices with all off-diagonal blocks equal to zero. Exploiting these
structures substantially accelerates the computation.


\section{Merging two DtN maps}
\label{sec5:merge}

Let $\tau$ denote a box in the tree with children $\alpha$ and $\beta$. In this
section, we demonstrate that if the DtN matrices $\mtx{T}^{\alpha}$ and $\mtx{T}^{\beta}$
for the children are known, then the DtN matrix $\mtx{T}^{\tau}$ can be constructed
via a purely local computation which we refer to as a ``merge'' operation.

\begin{figure}
\setlength{\unitlength}{1mm}
\begin{picture}(95,55)
\put(-15,00){\includegraphics[height=55mm]{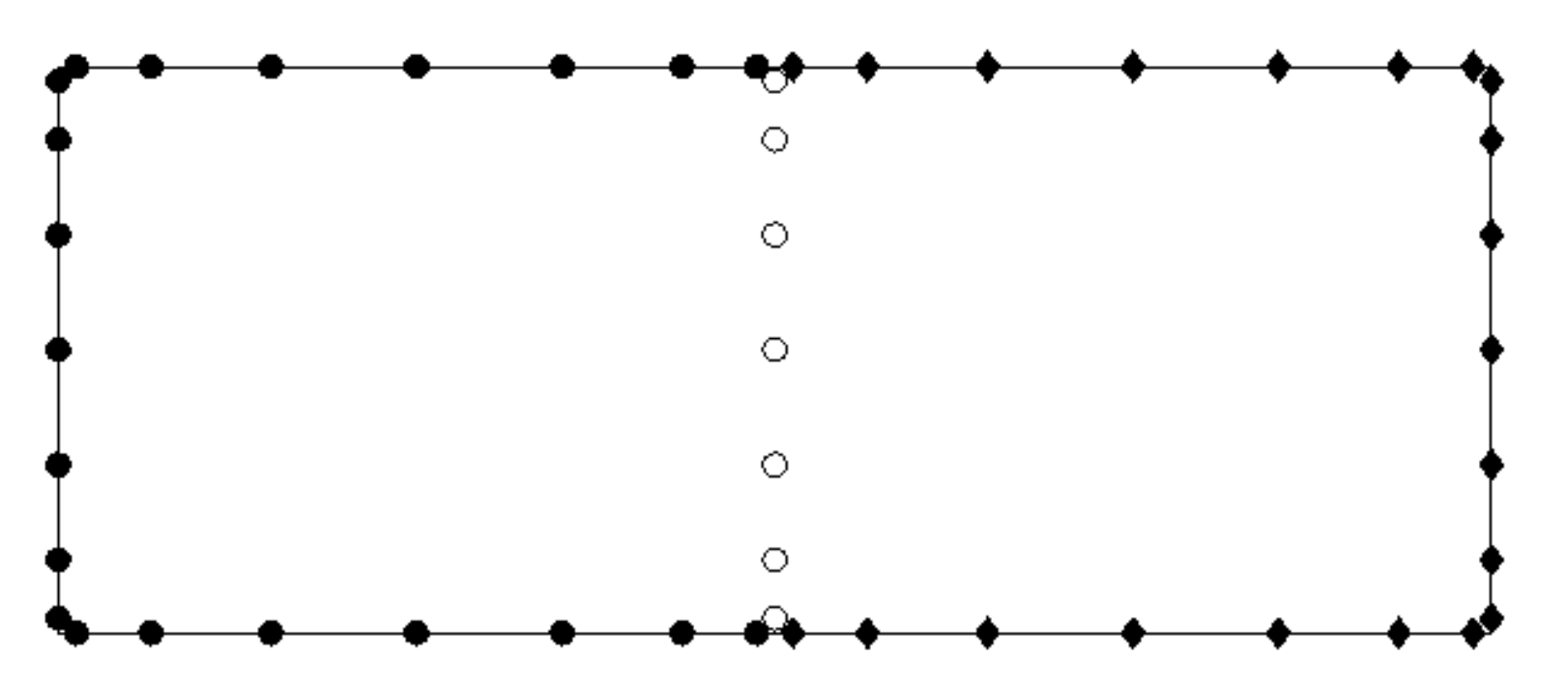}}
\put(18,25){$\Omega_{\alpha}$}
\put(74,25){$\Omega_{\beta}$}
\put(-9,25){$J_{1}$}
\put(100,25){$J_{2}$}
\put(42,25){$J_{3}$}
\end{picture}
\caption{Notation for the merge operation described in Section \ref{sec5:merge}.
The rectangular domain $\Omega$ is formed by two squares $\Omega_{\alpha}$ and $\Omega_{\beta}$.
The sets $J_{1}$ and $J_{2}$ form the exterior nodes (black), while
$J_{3}$ consists of the interior nodes (white).}
\label{fig5:siblings_notation}
\end{figure}
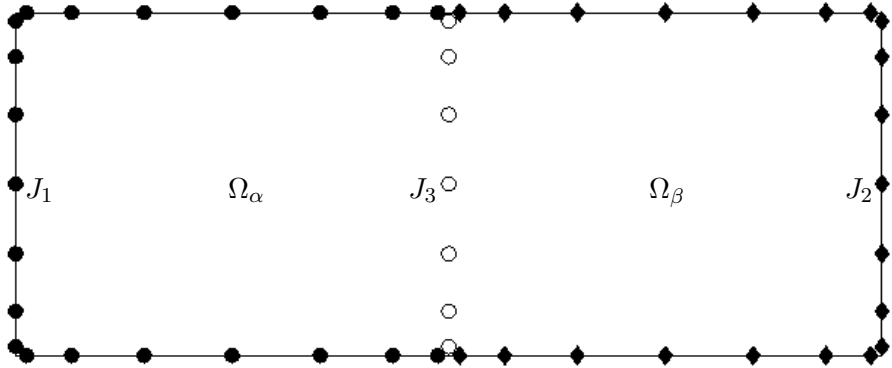

We start by introducing some notation:
Let $\Omega_{\tau}$ denote a box with children $\Omega_{\alpha}$ and
$\Omega_{\beta}$. For concreteness, let us assume that $\Omega_{\alpha}$ and
$\Omega_{\beta}$ share a vertical edge as shown in Figure \ref{fig5:siblings_notation},
so that
$$
\Omega_{\tau} = \Omega_{\alpha} \cup \Omega_{\beta}.
$$
We partition the points on $\partial\Omega_{\alpha}$ and $\partial\Omega_{\beta}$ into three sets:
\begin{tabbing}
\mbox{}\hspace{5mm}\= $J_{1}$ \hspace{4mm} \=
Boundary nodes of $\Omega_{\alpha}$ that are not boundary nodes of $\Omega_{\beta}$.\\
\> $J_{2}$ \> Boundary nodes of $\Omega_{\beta}$ that are not boundary nodes of $\Omega_{\alpha}$.\\
\> $J_{3}$ \> Boundary nodes of both $\Omega_{\alpha}$ and $\Omega_{\beta}$ that are \textit{not} boundary nodes of the union box
$\Omega_{\tau}$.
\end{tabbing}
Figure \ref{fig5:siblings_notation} illustrates the definitions of the $J_{k}$'s.
Let $u$ denote a solution to (\ref{eq5:basic}), with tabulated potential values $\uu$
and boundary fluxes $\vv$, as described in Section \ref{sec5:discretization}. Set
\begin{equation}
\label{eq5:defuiue}
\uu_{\rm i} = \uu_{3},
\qquad\mbox{and}\qquad
\uu_{\rm e} = \left[\begin{array}{c} \uu_{1} \\ \uu_{2} \end{array}\right].
\end{equation}
Recall that $\TT^{\alpha}$ and $\TT^{\beta}$ denote the
operators that map values of the potential $u$ on the boundary to values of
$\partial_{n}u$ on the boundaries
of the boxes $\Omega_{\alpha}$ and $\Omega_{\beta}$,
as described in Section \ref{sec5:leaf}.
The operators can be partitioned according to the numbering of nodes
in Figure \ref{fig5:siblings_notation}, resulting in the equations
\begin{equation}
\label{eq5:bittersweet1}
\left[\begin{array}{c}
\vv_{1}\\ \vv_{3}
\end{array}\right] =
\left[\begin{array}{ccc}
\mtx{T}_{1,1}^{\alpha} & \mtx{T}_{1,3}^{\alpha} \\
\mtx{T}_{3,1}^{\alpha} & \mtx{T}_{3,3}^{\alpha}
\end{array}\right]\,
\left[\begin{array}{c}
\uu_{1}\\ \uu_{3}
\end{array}\right],
\qquad {\rm and} \qquad
\left[\begin{array}{c}
\vv_{2}\\ \vv_{3}
\end{array}\right] =
\left[\begin{array}{ccc}
\mtx{T}_{2,2}^{\beta} & \mtx{T}_{2,3}^{\beta} \\
\mtx{T}_{3,2}^{\beta} & \mtx{T}_{3,3}^{\beta}
\end{array}\right]\,
\left[\begin{array}{c}
\uu_{2}\\ \uu_{3}
\end{array}\right].
\end{equation}

Our objective is now to construct a solution operator $\mtx{S}^{\tau}$
and a DtN matrix $\mtx{T}^{\tau}$ such that
\begin{align}
\label{eq5:desiredS}
\uu_{3} =&\ \mtx{S}^{\tau}\,
\left[\begin{array}{c}\uu_{1} \\ \uu_{2}\end{array}\right]\\
\label{eq5:desiredT}
\left[\begin{array}{c}\vv_{1} \\ \vv_{2}\end{array}\right]
=&\
\mtx{T}^{\tau}\,
\left[\begin{array}{c}\uu_{1} \\ \uu_{2}\end{array}\right].
\end{align}
To this end, we eliminate $\vct{v}_{3}$ from (\ref{eq5:bittersweet1}) and write the result as a single equation:
\begin{equation}
\label{eq5:premerge}
\left[\begin{array}{cc|c}
\TT^{\alpha}_{1,1} & \mtx{0} & \TT^{\alpha}_{1,3} \\
\mtx{0} & \TT^{\beta }_{2,2} & \TT^{\beta }_{2,3} \\ \hline
\TT^{\alpha}_{3,1} & -\TT^{\beta}_{3,2} & \TT^{\alpha}_{3,3} - \TT^{\beta}_{3,3}
\end{array}\right]\,
\left[\begin{array}{c}
\uu_{1} \\ \uu_{2} \\ \hline \uu_{3}
\end{array}\right]
=
\left[\begin{array}{c}
\vv_{1} \\ \vv_{2} \\ \hline \vct{0}
\end{array}\right],
\end{equation}
The last equation directly tells us that (\ref{eq5:desiredS}) holds with
\begin{equation}
\label{eq5:U_parent}
\mtx{S}^{\tau} =
\bigl(\TT^{\alpha}_{3,3} - \TT^{\beta}_{3,3}\bigr)^{-1}
\bigl[-\TT^{\alpha}_{3,1}\ \big|\ \TT^{\beta}_{3,2}].
\end{equation}
By eliminating $\uu_{3}$ from (\ref{eq5:premerge}) by forming a Schur complement, we also find that
(\ref{eq5:desiredT}) holds with
\begin{equation}
\label{eq5:merge}
\mtx{T}^{\tau} =
\left[\begin{array}{ccc}
\TT_{1,1}^{\alpha} & \mtx{0} \\
\mtx{0} & \TT_{2,2}^{\beta }
\end{array}\right] +
\left[\begin{array}{c}
\TT_{1,3}^{\alpha} \\
\TT_{2,3}^{\beta}
\end{array}\right]\,
\bigl(\TT^{\alpha}_{3,3} - \TT^{\beta}_{3,3}\bigr)^{-1}
\bigl[-\TT^{\alpha}_{3,1}\ \big|\ \TT^{\beta}_{3,2}\bigr].
\end{equation}


\section{The full hierarchical scheme}
\label{sec5:hierarchy}

At this point, we know how to construct the DtN operator for a leaf (Section
\ref{sec5:leaf}), and how to merge two such operators of neighboring
patches to form the DtN operator of their union (Section
\ref{sec5:merge}). We are ready to describe the full hierarchical
scheme for solving the Dirichlet problem (\ref{eq5:basic}). This scheme
takes the Dirichlet boundary data $f$, and constructs an approximation
to the solution $u$. The output is a vector $\uu$ that tabulates
approximations to $u$ at the Gaussian nodes $\{\pxx_{k}\}_{k=1}^{N}$
on all interior edges that were defined in Section \ref{sec5:discretization}.
To find $u$ at an arbitrary set of target points in $\Omega$, a post-processing
step described in Section \ref{sec5:postproc} can be used.

\subsection{The algorithm}
\label{sec5:thealgorithm}
Partition the domain into a hierarchical tree as described in Section \ref{sec5:discretization}.
Then execute a ``build stage'' in which we construct for each box $\tau$ the following two matrices:

\vspace{2mm}

\begin{itemize}
\item[$\mtx{S}^{\tau}$]
For a parent box $\tau$, $\mtx{S}^{\tau}$ is a solution operator that maps values
of $u$ on $\partial \Omega_{\tau}$ to values of $u$ at the interior nodes. In other
words, $\uu(I_{\rm i}^{\tau}) = \mtx{S}^{\tau}\,\uu(I_{\rm e}^{\tau})$.
(For a leaf $\tau$, $\mtx{S}^{\tau}$ is not defined.)

\vspace{2mm}

\item[$\TT^{\tau}$] The matrix that maps $\uu(I_{\rm e}^{\tau})$ (tabulating values of $u$ on
$\partial \Omega_{\tau}$) to $\vv(I_{\rm e}^{\tau})$ (tabulating values of $du/dn$).
In other words, $\vv(I_{\rm e}^{\tau}) = \TT^{\tau}\,\uu(I_{\rm e}^{\tau})$.
\end{itemize}

\vspace{2mm}

\noindent
(Recall that the index vectors $I_{\rm e}^{\tau}$ and $I_{\rm i}^{\tau}$
were defined in Section \ref{sec5:discretization}.)
The build stage consists of a single sweep over all nodes in the tree.
Any bottom-up ordering in which any parent box is processed after its
children can be used. For each leaf box $\tau$, an approximation
to the local DtN map $\mtx{T}^{\tau}$ is constructed using the procedure
described in Section \ref{sec5:leaf}. For a parent box
$\tau$ with children $\alpha$ and $\beta$, the matrices
$\mtx{S}^{\tau}$ and $\TT^{\tau}$ are formed from the DtN operators $\TT^{\alpha}$
and $\TT^{\beta}$ via the process described in Section \ref{sec5:merge}.
Algorithm 1 summarizes the build stage.


Once all the matrices $\{\mtx{S}^{\tau}\}_{\tau}$ have been formed, a vector $\uu$ holding approximations
to the solution $u$ of (\ref{eq5:basic}) can be constructed for all discretization points by
starting at the root box $\Omega$ and moving down the tree toward the leaf boxes.
The values of $\uu$ for the points on the boundary of $\Omega$ can be obtained by tabulating the
boundary function $f$.  When any box $\tau$ is processed, the value of $\uu$ is known for all nodes
on its boundary (i.e.~those listed in $I_{\rm e}^{\tau}$). The matrix $\mtx{S}^{\tau}$ directly
maps these values to the values of $\uu$ on the nodes in the interior of $\tau$ (i.e.~those
listed in $I_{\rm i}^{\tau}$). When all nodes have been processed, approximations to $u$ have
constructed for all tabulation nodes on interior edges.
Algorithm 2 summarizes the solve stage.


\begin{figure}
\fbox{
\begin{minipage}{140mm}
\begin{center}
\textsc{Algorithm 1} (build solution operators)
\end{center}

This algorithm builds the global Dirichlet-to-Neumann operator for (\ref{eq5:basic}).\\
It also builds all matrices $\mtx{S}^{\tau}$ required for constructing $u$ at any interior point.\\
It is assumed that if node $\tau$ is a parent of node $\sigma$, then $\tau < \sigma$.

\rule{\textwidth}{0.5pt}

\begin{tabbing}
\mbox{}\hspace{7mm} \= \mbox{}\hspace{6mm} \= \mbox{}\hspace{6mm} \= \mbox{}\hspace{6mm} \= \mbox{}\hspace{6mm} \= \kill
(1)\> \textbf{for} $\tau = N_{\rm boxes},\,N_{\rm boxes}-1,\,N_{\rm boxes}-2,\,\dots,\,1$\\
(2)\> \> \textbf{if} ($\tau$ is a leaf)\\
(3)\> \> \> Construct $\mtx{T}^{\tau}$ via the process described in Section \ref{sec5:leaf}.\\
(4)\> \> \textbf{else}\\
(5)\> \> \> Let $\alpha$ and $\beta$ be the children of $\tau$.\\
(6)\> \> \> Split $I_{\rm e}^{\alpha}$ and $I_{\rm e}^{\beta}$ into vectors $I_{1}$, $I_{2}$, and $I_{3}$ as shown in Figure \ref{fig5:siblings_notation}.\\
(7)\> \> \> $\mtx{S}^{\tau} = \bigl(\TT^{\alpha}_{3,3} - \TT^{\beta}_{3,3}\bigr)^{-1}
                           \bigl[-\TT^{\alpha}_{3,1}\  \big|\
                                  \TT^{\beta}_{3,2}\bigr]$\\
(8)\> \> \> $\TT^{\tau} = \left[\begin{array}{ccc}
                          \mtx{T}_{1,1}^{\alpha} & \mtx{0}\\
                          \mtx{0} & \mtx{T}_{2,2}^{\beta }
                          \end{array}\right] +
                    \left[\begin{array}{c}
                          \TT_{1,3}^{\alpha} \\
                          \TT_{2,3}^{\beta}
                          \end{array}\right]\,\mtx{S}^{\tau}$.\\
(9)\> \> \> Delete $\TT^{\alpha}$ and $\TT^{\beta}$.\\
(10)\> \> \textbf{end if}\\
(11)\> \textbf{end for}
\end{tabbing}
\end{minipage}}
\end{figure}

\begin{figure}
\fbox{
\begin{minipage}{140mm}
\begin{center}
\textsc{Algorithm 2} (solve BVP once solution operator has been built)
\end{center}

This program constructs an approximation $\uu$ to the solution $u$ of (\ref{eq5:basic}).\\
It assumes that all matrices $\mtx{S}^{\tau}$ have already been constructed in a pre-computation.\\
It is assumed that if node $\tau$ is a parent of node $\sigma$, then $\tau < \sigma$.

\rule{\textwidth}{0.5pt}

\begin{tabbing}
\mbox{}\hspace{7mm} \= \mbox{}\hspace{6mm} \= \mbox{}\hspace{6mm} \= \mbox{}\hspace{6mm} \= \mbox{}\hspace{6mm} \= \kill
(1)\> $\uu(k) = f(\pxx_{k})$ for all $k \in I_{\rm e}^{1}$.\\
(2)\> \textbf{for} $\tau = 1,\,2,\,3,\,\dots,\,N_{\rm boxes}$\\
(3)\> \> \textbf{if} ($\tau$ is a parent)\\
(4)\> \> \> $\uu(I_{\rm i}^{\tau}) = \mtx{S}^{\tau}\,\uu(I_{\rm e}^{\tau})$.\\
(5)\> \> \textbf{end if}\\
(6)\> \textbf{end for}
\end{tabbing}

\vspace{1mm}

\textit{Remark: This algorithm outputs the solution on the Gaussian nodes on box boundaries.
To get the solution at other points, use the method described in Section \ref{sec5:postproc}.}

\end{minipage}}
\end{figure}

\begin{remark}
The merge stage is exact when performed in exact arithmetic. The only
approximation involved is the approximation of the solution $u$ on a leaf by its
interpolating polynomial.
\end{remark}

\begin{remark}
\label{remark:refine}
To keep the presentation simple, we consider in this note only the case
of a uniform computational grid. Such grids are obviously not well suited
to situations where the regularity of the solution changes across the domain.
The method described can \textit{in principle} be modified to handle locally refined
grids quite easily. A complication is that the tabulation nodes for
two touching boxes will typically not coincide, which requires the introduction of
specialized interpolation operators. Efficient refinement strategies also
require the development of error indicators that identify the regions where
the grid need to be refined. This is work in progress, and will be reported
at a later date. We observe that our introduction of Gaussian nodes on the
internal boundaries (as opposed to the Chebyshev nodes used in
\cite{2012_spectralcomposite}) makes re-interpolation much
easier.
\end{remark}

\subsection{Asymptotic complexity}
\label{sec5:complexity}
In this section, we determine the asymptotic complexity of the direct solver.
Let $N_{\rm leaf} = 4\ngauss$ denote the number of Gaussian nodes on the boundary
of a leaf box, and let $q^2$ denote the number of Chebychev nodes used in the leaf computation.
Let $L$ denote the number of levels in the binary tree.  This means there are
$4^L$ boxes.  Thus the total number of discretization nodes $N$ is approximately $4^L\ngauss = \frac{(2^L\ngauss)^2}{\ngauss}$.
(To be exact, $N = 2^{2L+1}\ngauss+2^{L+1}\ngauss$.)

The cost to process one leaf is approximately $O(q^6)$.  Since there are $\frac{N}{q^2}$ leaf boxes, the total
cost of pre-computing approximate DtN operators for all the bottom level is
$\frac{N}{q^2}\times q^6\sim N\ngauss^4$.

Next, consider the cost of constructing the DtN map on level $\ell$ via
the merge operation described in Section \ref{sec5:merge}.  For each box
on the level $\ell$, the operators $\TT^\tau$ and $\mtx{S}^\tau$ are constructed
via (\ref{eq5:U_parent}) and (\ref{eq5:U_parent}).  These operations involve
matrices of size roughly $2^{-\ell}N^{0.5}\times  2^{-\ell}N^{0.5}$.  Since there are
$4^\ell$ boxes per level.  The cost on level $\ell$ of the merge is
$$4^\ell\times \left(2^{-\ell}N^{0.5}\right)^3 \sim 2^{-\ell}N^{1.5}.$$

The total cost for all the merge procedures has complexity
$$\sum_{\ell=1}^L 2^{-\ell}N^{1.5} \sim N^{1.5}.$$

Finally, consider the cost of the downwards sweep which solves for the
interior unknowns.  For any non-leaf box $\tau$ on level $\ell$, the size of $\mtx{S}^\tau$ is
$2^l\ngauss \times  2^l(6\ngauss)$ which is approximately $\sim 2^{-\ell}N^{0.5} \times 2^{-\ell}N^{0.5}$.
 Thus the cost of applying $\mtx{S}^\tau$ is roughly $(2^{-\ell}N^{0.5})^2 = 2^{-2\ell}N$.
So the total cost of the solve step has complexity
$$ \sum_{l=0}^{L-1}2^{2\ell}2^{-2\ell}N \sim N\log N.$$

In \cite{2013_martinsson_DtN_linearcomplexity}, we explain how to exploit structure in
the matrices $\TT$ and $\mtx{S}$ to improve the computational
cost of both the precomputation and the solve steps.

\subsection{Post-processing}
\label{sec5:postproc}

The direct solver in Algorithm 1 constructs approximations
to the solution $u$ of (\ref{eq5:basic}) at tabulation nodes at all interior edges.
Once these are available, it is easy to construct an approximation to $u$ at
an arbitrary point. To illustrate the process, suppose that we seek an approximation
to $u(\pyy)$, where $\pyy$ is a point located in a leaf $\tau$. We have values of $u$
tabulated at Gaussian nodes on $\partial \Omega_{\tau}$. These can easily be re-interpolated
to the Chebyshev nodes on $\partial \Omega_{\tau}$. Then $u$ can be reconstructed at
the interior Chebyshev nodes via the formula (\ref{eq5:leaf_solve}); observe that
the local solution operator $-\mtx{A}_{\rm i,i}^{-1}\mtx{A}_{\rm i,e}$ was built when
the leaf was originally processed and can be simply retrieved from memory (assuming
enough memory is available). Once $u$ is tabulated at the Chebyshev grid on $\Omega_{\tau}$, it is
trivial to interpolate it to $\pyy$ or any other point.

\section{Body loads}
\label{sec5:body}

Now that we have described how to solve our basic boundary value problem
\begin{equation}
\left\{\begin{aligned}
\mbox{}[Au](\pxx) =&\ g(\pxx),\qquad &\pxx \in \Omega,\\
   u(\pxx) =&\ f(\pxx),\qquad &\pxx \in \Gamma,
\end{aligned}\right.
\end{equation}
for the special case where $g=0$, we will next consider the more general
case that includes a body load. Only minor modifications are required to
the basic scheme, but note that the resulting method requires substantially
more memory. The techniques presented here were developed jointly with Tracy
Babb of CU-Boulder, cf.~\cite{2014_haut_hyperbolic}.

\subsection{Notation}
When handling body loads, we will extensively switch between the Chebyshev
and Gaussian grids, so we need to introduce some additional notation.

Let $\{\pyy_{j}\}_{j=1}^{M}$ denote the global grid obtained by putting
down a $\ncheb \times \ncheb$ tensor product grid of Chebyshev nodes on
each leaf. For a leaf $\tau$, let $I_{\rm c}^{\tau}$ denote an index vector
pointing to the nodes in $\{\pyy_{j}\}_{j=1}^{M}$ that lie on leaf $\tau$.
We partition this index vector into exterior and interior nodes as follows
$$
I_{\rm c}^{\tau} = I_{\rm ce}^{\tau} \cup I_{\rm ci}^{\tau}.
$$
To avoid confusion with the index vectors pointing into the grid of
Gaussian boundary functions, we rename these index vectors as follows:
\begin{description}
\item[$I_{\rm ge}^{\tau}$] An index vector marking the (Gauss) nodes in
$\{\pxx_{i}\}_{i=1}^{N}$ that lie on $\partial\Omega_{\tau}$.
\item[$I_{\rm gi}^{\tau}$] For a parent node $\tau$, this is an index vector
marking the (Gauss) nodes in $\{\pxx_{i}\}_{i=1}^{N}$ that lie on the
``interior'' boundary of $\tau$. For a leaf node, this vector is not defined.
\end{description}

The fact that we use two spectral grids also leads to a need to distinguish
between different vectors tabulating approximate values. We have:

\begin{tabular}{l|l|l}
                          & $\tau$ is a leaf          & $\tau$ is a parent \\ \hline
$\vct{u}_{\rm c }^{\tau}$ & $u$ tabulated on Chebyshev nodes          & --- \\
$\vct{u}_{\rm ce}^{\tau}$ & $u$ tabulated on Chebyshev exterior nodes & --- \\
$\vct{u}_{\rm ci}^{\tau}$ & $u$ tabulated on Chebyshev interior nodes & --- \\
$\vct{u}_{\rm ge}^{\tau}$ & $u$ tabulated on Gaussian  exterior nodes & $u$ tabulated on Gaussian exterior nodes \\
$\vct{u}_{\rm gi}^{\tau}$ & ---                                       & $u$ tabulated on Gaussian interior nodes
\end{tabular}

\subsection{Leaf computation}
Let $\tau$ be a leaf. We split the solution $u$ to the equation
\begin{equation}
\label{eq5:bodyload}
\left\{\begin{aligned}
Au(\pxx) =&\ g(\pxx),\qquad &\pxx \in \Omega_{\tau},\\
 u(\pxx) =&\ \psi(\pxx),\qquad &\pxx \in \Gamma_{\tau},
\end{aligned}\right.
\end{equation}
as
$$
u = w + \phi
$$
where $w$ is a \textit{particular solution}
\begin{equation}
\label{eq5:part}
\left\{\begin{aligned}
Aw(\pxx) =&\ g(\pxx),\qquad &\pxx \in \Omega_{\tau},\\
 w(\pxx) =&\ 0,\qquad &\pxx \in \Gamma_{\tau},
\end{aligned}\right.
\end{equation}
and where $\phi$ is a \textit{homogeneous solution}
\begin{equation}
\label{eq5:hom}
\left\{\begin{aligned}
A\phi(\pxx) =&\ 0,\qquad &\pxx \in \Omega_{\tau},\\
 \phi(\pxx) =&\ \psi(\pxx),\qquad &\pxx \in \Gamma_{\tau}.
\end{aligned}\right.
\end{equation}
We can now write the Neumann data for $u$ as
$$
\partial_{n}u|_{\Gamma_{\tau}}
=
\partial_{n}w|_{\Gamma_{\tau}} +
\partial_{n}\phi|_{\Gamma_{\tau}}
=
\partial_{n}w|_{\Gamma_{\tau}} +
T\phi|_{\Gamma_{\tau}}
=
\partial_{n}w|_{\Gamma_{\tau}} +
T\,\psi,
$$
where, as before, $T$ is the NfD operator. Our objective is therefore to find a matrix
that maps the given body load $g$ on $\Omega_{\tau}$ to the Neumann data of $w$. We do this
calculation on the Chebyshev grid. Discretizing (\ref{eq5:part}), and collocating on
the internal nodes, we find
$$
\mtx{A}_{\rm ci,ce}\vct{w}_{\rm ce} + \mtx{A}_{\rm ci,ci}\vct{w}_{\rm ci} = \vct{g}_{\rm ci}.
$$
But now observe that $\vct{w}_{\rm ce} = 0$, so the particular solution is given by
$$
\vct{w}_{\rm c} =
\left[\begin{array}{c} \vct{w}_{\rm ce} \\ \vct{w}_{\rm ci} \end{array}\right] =
\mtx{F}_{\rm c,ci}\vct{g}_{\rm ci},
\qquad\mbox{where}\qquad
\mtx{F}_{\rm c,ci} =
\left[\begin{array}{c} \mtx{0} \\ \mtx{A}_{\rm ci,ci}^{-1} \end{array}\right].
$$
Let $\vct{h}_{\rm ge}$ denote the Neumann data for $\vct{w}$, tabulated at the Gaussian
nodes on the boundary. It follows that
$$
\vct{h}_{\rm ge} =
\underbrace{\mtx{D}_{\rm ge,c}\mtx{F}_{\rm c,ci}}_{=:\mtx{H}_{\rm ge,ci}}\vct{g}_{\rm ci},
$$
where $\mtx{D}_{\rm ge,c}$ is the operator that differentiates from the full Chebyshev
grid, and then interpolates to the Gaussian exterior nodes. (Using the notation from
Section \ref{sec5:merge}, $\mtx{D}_{\rm ge,c} = \mtx{L}_{4}\mtx{L}_{3}$.)

Once the boundary data of $u$ on $\Omega_{\tau}$ is given, the total solution on the
Chebyshev grid is given by
$$
\vct{u}_{\rm c} =
\mtx{S}_{\rm c,ge}\vct{u}_{\rm ge} + \mtx{F}_{\rm c,ci}\vct{g}_{\rm ci} =
\mtx{S}_{\rm c,ge}\vct{u}_{\rm ge} + \vct{w}_{\rm c},
$$
where $\mtx{S}_{\rm c,ge}$ is the solution operator given by
$$
\mtx{S}_{\rm c,ge} =
\left[\begin{array}{c}\mtx{I}\\ -\mtx{A}_{\rm ci,ci}^{-1}\mtx{A}_{\rm ci,ce}\end{array}\right]
\mtx{I}_{\rm ce,ge},
$$
where, in turn, $\mtx{I}_{\rm ce,ge}$ is the interpolation operator that maps from the
boundary Gauss nodes to the boundary Chebyshev nodes.

\subsection{Merge operator}
Let $\tau$ be a parent node with children $\alpha$ and $\beta$. The local equilibrium equations read
\begin{align}
\label{eq5:bitter1}
\left[\begin{array}{c}
\vv_{1}\\ \vv_{3}
\end{array}\right]
=&\
\left[\begin{array}{ccc}
\mtx{T}_{1,1}^{\alpha} & \mtx{T}_{1,3}^{\alpha} \\
\mtx{T}_{3,1}^{\alpha} & \mtx{T}_{3,3}^{\alpha}
\end{array}\right]\,
\left[\begin{array}{c}
\uu_{1}\\ \uu_{3}
\end{array}\right]
+
\left[\begin{array}{c} \vct{h}_{1}^{\alpha} \\ \vct{h}_{3}^{\alpha} \end{array}\right],
\\
\left[\begin{array}{c}
\vv_{2}\\ \vv_{3}
\end{array}\right] =&\
\left[\begin{array}{ccc}
\mtx{T}_{2,2}^{\beta} & \mtx{T}_{2,3}^{\beta} \\
\mtx{T}_{3,2}^{\beta} & \mtx{T}_{3,3}^{\beta}
\end{array}\right]\,
\left[\begin{array}{c}
\uu_{2}\\ \uu_{3}
\end{array}\right]
+
\left[\begin{array}{c} \vct{h}_{2}^{\beta} \\ \vct{h}_{3}^{\beta} \end{array}\right].
\end{align}
Combine the two equations for $\vv_{3}$ to obtain the equation
$$
\mtx{T}_{3,1}^{\alpha}\,\uu_{1} + \mtx{T}_{3,3}^{\alpha}\,\uu_{3} + \vct{h}_{3}^{\alpha}
=
\mtx{T}_{3,2}^{\beta }\,\uu_{2} + \mtx{T}_{3,3}^{\beta }\,\uu_{3} + \vct{h}_{3}^{\beta}.
$$
This gives [check for sign errors!]
\begin{equation}
\label{eq5:swiss}
\uu_{3} =
\bigl(\TT^{\alpha}_{3,3} - \TT^{\beta}_{3,3}\bigr)^{-1}
\bigl( \TT^{\beta }_{3,2}\uu_{2}
      -\TT^{\alpha}_{3,1}\uu_{1}
      +\vct{h}_{3}^{\beta}
      -\vct{h}_{3}^{\alpha}
      \bigr)
\end{equation}
Inserting (\ref{eq5:swiss}) back into (\ref{eq5:bitter1}) we find
\begin{align*}
\left[\begin{array}{c} \vv_{1} \\ \vv_{2} \end{array}\right] =&\
\left(\left[\begin{array}{ccc}
\TT_{1,1}^{\alpha} & \mtx{0} \\
\mtx{0} & \TT_{2,2}^{\beta }
\end{array}\right] +
\left[\begin{array}{c}
\TT_{1,3}^{\alpha} \\
\TT_{2,3}^{\beta}
\end{array}\right]\,
\bigl(\TT^{\alpha}_{3,3} - \TT^{\beta}_{3,3}\bigr)^{-1}
\bigl[-\TT^{\alpha}_{3,1}\ \big|\ \TT^{\beta}_{3,2}].
\right)
\left[\begin{array}{c} \uu_{1} \\ \uu_{2} \end{array}\right] +\\
&\
\left[\begin{array}{c} \vct{h}^{\alpha}_{1} \\ \vct{h}^{\beta}_{2} \end{array}\right] +
\left[\begin{array}{c}
\TT_{1,3}^{\alpha} \\
\TT_{2,3}^{\beta}
\end{array}\right]\,
\bigl(\TT^{\alpha}_{3,3} - \TT^{\beta}_{3,3}\bigr)^{-1}
\bigl(\vct{h}_{3}^{\beta}-\vct{h}_{3}^{\alpha}\bigr).
\end{align*}
We now define operators
\begin{align*}
\mtx{X}^{\tau}_{\rm gi,gi} =&\ \bigl(\TT^{\alpha}_{3,3} - \TT^{\beta}_{3,3}\bigr)^{-1},\\
\mtx{S}^{\tau}_{\rm gi,ge} =&\
\bigl(\TT^{\alpha}_{3,3} - \TT^{\beta}_{3,3}\bigr)^{-1}\bigl[-\TT^{\alpha}_{3,1}\ \big|\ \TT^{\beta}_{3,2}] =
\mtx{X}^{\tau}_{\rm gi,gi}\bigl[-\TT^{\alpha}_{3,1}\ \big|\ \TT^{\beta}_{3,2}],\\
\mtx{T}^{\tau}_{\rm ge,ge} =&
\left[\begin{array}{ccc}
\TT_{1,1}^{\alpha} & \mtx{0} \\
\mtx{0} & \TT_{2,2}^{\beta }
\end{array}\right]
\bigl(\TT^{\alpha}_{3,3} - \TT^{\beta}_{3,3}\bigr)^{-1}\bigl[-\TT^{\alpha}_{3,1}\ \big|\ \TT^{\beta}_{3,2}]
=
\left[\begin{array}{ccc}
\TT_{1,1}^{\alpha} & \mtx{0} \\
\mtx{0} & \TT_{2,2}^{\beta }
\end{array}\right]\mtx{S}^{\tau}_{\rm gi,ge}.
\end{align*}
Then in the upwards pass in the solve, we compute
\begin{align*}
\vct{w}_{\rm gi}^{\tau} =&\ \mtx{X}^{\tau}\bigl(\vct{h}_{3}^{\beta}-\vct{h}_{3}^{\alpha}\bigr)\\
\vct{h}_{\rm ge}^{\tau} =&\ \left[\begin{array}{c}
\TT_{1,3}^{\alpha} \\
\TT_{2,3}^{\beta}
\end{array}\right]\vct{w}_{\rm gi}^{\tau},
\end{align*}
and in the solve stage, we recover $\vct{u}_{\rm gi}^{\tau}$ via
$$
\vct{u}^{\tau}_{\rm gi} = \mtx{S}_{\rm gi,ge}\vct{u}^{\tau}_{\rm ge} + \vct{w}_{\rm gi}^{\tau}.
$$

\begin{remark}[Physical interpretation of merge]
The quantities $\vct{w}^{\tau}_{\rm gi}$ and $\vct{h}_{\rm ge}^{\tau}$ have a simple
physical meaning. The vector $\vct{w}_{\rm gi}^{\tau}$ introduced above is simply a
tabulation of the particular solution $w^{\tau}$ associated with $\tau$ on the interior
boundary $\Gamma_{3}$, and $\vct{h}_{\rm ge}^{\tau}$ is the normal derivative of $w^{\tau}$.
To be precise, $w^{\tau}$ is the solution to the inhomogeneous problem, cf.~(\ref{eq5:part})
\begin{equation}
\label{eq5:part_tau}
\left\{\begin{aligned}
Aw^{\tau}(\pxx) =&\ g(\pxx),\qquad &\pxx \in \Omega_{\tau},\\
 w^{\tau}(\pxx) =&\ 0,\qquad &\pxx \in \Gamma_{\tau}.
\end{aligned}\right.
\end{equation}
We can re-derive the formula for $w|_{\Gamma_{3}}$ using the original mathematical operators
as follows: First observe that for $\pxx \in \Omega^{\alpha}$, we have
$A(w^{\tau} - w^{\alpha}) = g - g = 0$, so the NfD operator $T^{\alpha}$ applies to the
function $w^{\tau} - w^{\alpha}$:
$$
T^{\alpha}_{31}(w_{1}^{\tau} - w_{1}^{\alpha}) +
T^{\alpha}_{33}(w_{3}^{\tau} - w_{3}^{\alpha}) =
(\partial_{n}w^{\tau})|_{3} - (\partial_{n}w^{\alpha})|_{3}
$$
Use that $w_{1}^{\tau} = w_{1}^{\alpha} = w_{3}^{\alpha} = 0$, and that $(\partial_{n}w^{\alpha})|_{3} = h_{3}^{\alpha}$ to get
\begin{equation}
\label{eq5:legomine1}
T^{\alpha}_{33}w_{3}^{\tau} = (\partial_{n}w^{\tau})|_{3} - h_{3}^{\alpha}.
\end{equation}
Analogously, we get
\begin{equation}
\label{eq5:legomine2}
T^{\beta}_{33}w_{3}^{\tau} = (\partial_{n}w^{\tau})|_{3} - h_{3}^{\beta}.
\end{equation}
Combine (\ref{eq5:legomine1}) and (\ref{eq5:legomine2}) to eliminate $(\partial_{n}w^{\tau})|_{3}$ and obtain
$$
\bigl(T^{\alpha}_{33} - T^{\beta}_{33}\bigr) w_{3}^{\tau} =
 - h_{3}^{\alpha} + h_{3}^{\beta}.
$$
Observe that in effect, we can write the particular solution $w^{\tau}$ as
$$
w^{\tau}(\pxx) = \left\{\begin{array}{ll}
w^{\alpha}(\pxx) + \hat{w}^{\tau}(\pxx)\quad& \pxx \in \Omega^{\alpha},\\
w^{\beta }(\pxx) + \hat{w}^{\tau}(\pxx)\quad& \pxx \in \Omega^{\beta},
\end{array}\right.
$$
The function $w^{\tau}$ must of course be smooth across $\Gamma_{3}$, so the
function $\hat{w}^{\tau}$ must have a jump that exactly offsets the discrepancy
in the derivatives of $w^{\alpha}$ and $w^{\beta}$. This jump is precisely of
size $h^{\alpha} - h^{\beta}$.
\end{remark}


\begin{figure}
\fbox{
\begin{minipage}{140mm}
\begin{center}
\textsc{Algorithm 3} (Build stage for problems with body load)
\end{center}

\vspace{2mm}

\begin{minipage}{135mm}
This algorithms build all solution operators required to solve the non-homogeneous
BVP (\ref{eq5:bodyload}).
It is assumed that if node $\tau$ is a parent of node $\sigma$, then $\tau < \sigma$.
\end{minipage}

\rule{\textwidth}{0.5pt}

\begin{tabbing}
\mbox{}\hspace{6mm} \= \mbox{}\hspace{6mm} \= \mbox{}\hspace{6mm} \= \mbox{}\hspace{6mm} \= \kill
\textbf{for} $\tau = N_{\rm boxes},\,N_{\rm boxes}-1,\,N_{\rm boxes}-2,\,\dots,\,1$\\
\> \textbf{if} ($\tau$ is a leaf)\\[1mm]
\> \> $\mtx{F}_{\rm c,ci}^{\tau} = \left[\begin{array}{c}\mtx{0}\\ \mtx{A}_{\rm ci,ci}^{-1}\end{array}\right]$
\quad\textit{[pot.] $\leftarrow$ [body load]}\\[1mm]
\> \> $\mtx{H}_{\rm ge,ci}^{\tau} = \mtx{D}_{\rm ge,c}\mtx{F}_{\rm c,ci}^{\tau}$
\quad\textit{[deriv.] $\leftarrow$ [body load]}\\[1mm]
\> \> $\mtx{S}_{\rm c,ge}^{\tau} = \left[\begin{array}{c}\mtx{I}\\ -\mtx{A}_{\rm ci,ci}^{-1}\mtx{A}_{\rm ci,ce}\end{array}\right]\,
\mtx{I}_{\rm ce,ge}$
\quad\textit{[pot.] $\leftarrow$ [pot.]}\\[1mm]
\> \> $\mtx{T}_{\rm ge,ge}^{\tau} = \mtx{D}_{\rm ge,c}\mtx{S}_{\rm c,ge}$
\quad\textit{[deriv.] $\leftarrow$ [pot.] (NfD operator)}\\[1mm]
\> \textbf{else}\\
\> \> Let $\alpha$ and $\beta$ be the children of $\tau$.\\
\> \> Partition $I_{\rm e}^{\alpha}$ and $I_{\rm e}^{\beta}$ into vectors $I_{1}$, $I_{2}$, and $I_{3}$ as shown in Figure \ref{fig5:siblings_notation}.\\[1mm]
\> \> $\mtx{X}_{\rm gi,gi}^{\tau} = \bigl(\TT^{\alpha}_{3,3} - \TT^{\beta}_{3,3}\bigr)^{-1}$
\quad\textit{[pot.] $\leftarrow$ [deriv.]}\\[1mm]
\> \> $\mtx{S}_{\rm gi,ge}^{\tau} = \mtx{X}_{\rm gi,gi}^{\tau}
                           \bigl[-\TT^{\alpha}_{3,1}\  \big|\
                                  \TT^{\beta}_{3,2}\bigr]$
\quad\textit{[pot.] $\leftarrow$ [pot.]}\\[1mm]
\> \> $\TT_{\rm ge,ge}^{\tau} = \left[\begin{array}{ccc}
                          \TT_{1,1}^{\alpha} & \mtx{0}\\
                          \mtx{0} & \TT_{2,2}^{\beta }
                          \end{array}\right] +
                    \left[\begin{array}{c}
                          \TT_{1,3}^{\alpha} \\
                          \TT_{2,3}^{\beta}
                          \end{array}\right]\,\UU_{\rm gi,ge}^{\tau}$
\quad\textit{[deriv.] $\leftarrow$ [pot.] (NfD operator)}.\\[1mm]
\> \textbf{end if}\\
\textbf{end for}
\end{tabbing}
\end{minipage}}
\caption{Build stage.}
\label{fig5:precomp_with_bodyload}
\end{figure}


\begin{figure}
\fbox{
\begin{minipage}{140mm}
\begin{center}
\textsc{Algorithm 4} (Solver for problems with body load)
\end{center}

\vspace{2mm}

\begin{minipage}{135mm}
This program constructs an approximation $\uu$ to the solution $u$ of (\ref{eq5:basic}).
It assumes that all matrices required to represent the solution operator have already
been constructed using Algorithm 3.
It is assumed that if node $\tau$ is a parent of node $\sigma$, then $\tau < \sigma$.
\end{minipage}

\rule{\textwidth}{0.5pt}

\textit{\textbf{Upwards pass --- construct all particular solutions:}}
\begin{tabbing}
\mbox{}\hspace{6mm} \= \mbox{}\hspace{6mm} \= \mbox{}\hspace{6mm} \= \mbox{}\hspace{6mm} \= \kill
\textbf{for} $\tau = N_{\rm boxes},\,N_{\rm boxes}-1,\,N_{\rm boxes}-2,\,\dots,\,1$\\
\> \textbf{if} ($\tau$ is a leaf)\\
\> \> \textit{\# Compute the derivatives of the local particular solution.}\\
\> \> $\vct{h}_{\rm ge}^{\tau} = \mtx{H}_{\rm ge,ci}^{\tau}\,\vct{g}_{\rm ci}^{\tau}$\\
\> \textbf{else} \\
\> \> Let $\alpha$ and $\beta$ denote the children of $\tau$.\\
\> \> \textit{\# Compute the local particular solution.}\\
\> \> $\vct{w}_{\rm gi}^{\tau} = \mtx{X}_{\rm gi,gi}^{\tau}\,\bigl(-\vct{h}_{3}^{\alpha} + \vct{h}_{3}^{\beta}\bigr).$\\
\> \> \textit{\# Compute the derivatives of the local particular solution.}\\
\> \> $\vct{h}_{\rm ge}^{\tau} =
       \left[\begin{array}{c}\vct{h}_{1}^{\alpha} \\ \vct{h}_{2}^{\beta}\end{array}\right] +
       \left[\begin{array}{c}\mtx{T}_{1,3}^{\alpha} \\ \mtx{T}_{2,3}^{\beta}\end{array}\right]\,\vct{w}_{\rm gi}^{\tau}.$\\
\> \textbf{end if}\\
\textbf{end for}
\end{tabbing}

\lsp

\textit{\textbf{Downwards pass --- construct all potentials:}}
\begin{tabbing}
\mbox{}\hspace{6mm} \= \mbox{}\hspace{6mm} \= \mbox{}\hspace{6mm} \= \mbox{}\hspace{6mm} \= \kill
\textit{\# Use the provided Dirichlet data to set the solution on the exterior of the root.}\\
$\uu(k) = f(\pxx_{k})$ for all $k \in I_{\rm e}^{1}$.\\
\textbf{for} $\tau = 1,\,2,\,3,\,\dots,\,N_{\rm boxes}$\\
\> \textbf{if} ($\tau$ is a parent)\\
\> \> \textit{\# Add the homogeneous term and the particular term.}\\
\> \> $\uu_{\rm gi}^{\tau} = \UU_{\rm gi,ge}^{\tau}\,\uu_{\rm ge}^{\tau} + \vct{w}_{\rm gi}^{\tau}$.\\
\> \textbf{else}\\
\> \> \textit{\# Add the homogeneous term and the particular term.}\\
\> \> $\uu_{\rm c}^{\tau} = \UU_{\rm c,ge}^{\tau}\,\uu_{\rm ge}^{\tau} +
                            \mtx{F}_{\rm c,ci}^{\tau}\,\vct{g}_{\rm ci}^{\tau}$.\\
\> \textbf{end}\\
\textbf{end for}
\end{tabbing}
\end{minipage}}
\caption{Solve stage.}
\label{fig5:solver_with_bodyload}
\end{figure}

\clearpage

\begin{appendix}
\section{A graphical illustration of the algorithm}
\label{sec5:cartoon}

This section provides an illustrated overview of the hierarchical
merge process described in detail in Section \ref{sec5:thealgorithm}.
The figures illustrate a situation
in which a square domain $\Omega = [0,1]^{2}$ is split into $4 \times 4$
leaf boxes on the finest level, and a $6\times 6$ spectral grid is
used in each leaf.

\vspace{2mm}

\noindent
\textbf{Step 1:} Partition the box $\Omega$ into $16$ small boxes that
each holds an $6\times 6$ Cartesian mesh of Chebyshev nodes.
For each box, identify the internal nodes (marked in blue) and
eliminate them as described in Section \ref{sec5:leaf}. Construct
the solution operator $\UU$, and the DtN operators
encoded in the matrices $\VV$ and $\WW$.

\begin{center}
\setlength{\unitlength}{1mm}
\begin{picture}(130,57)
\put(00,00){\includegraphics[height=55mm]{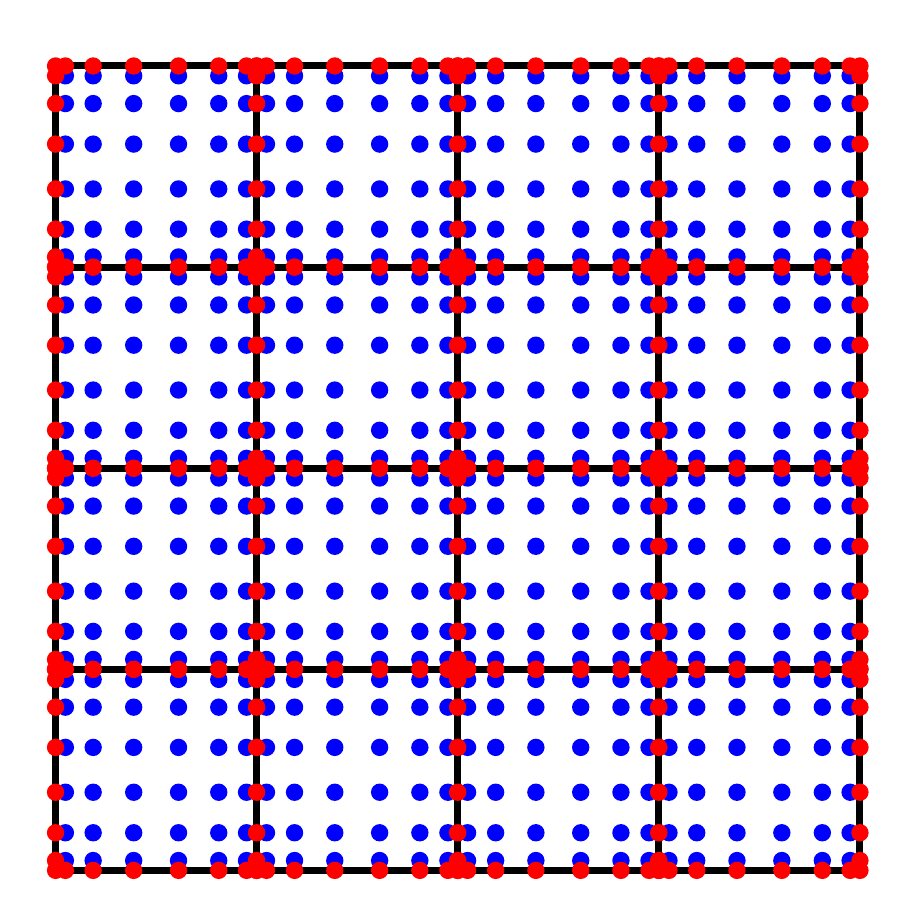}}
\put(63,27){$\Rightarrow$}
\put(60,32){Step 1}
\put(73,00){\includegraphics[height=55mm]{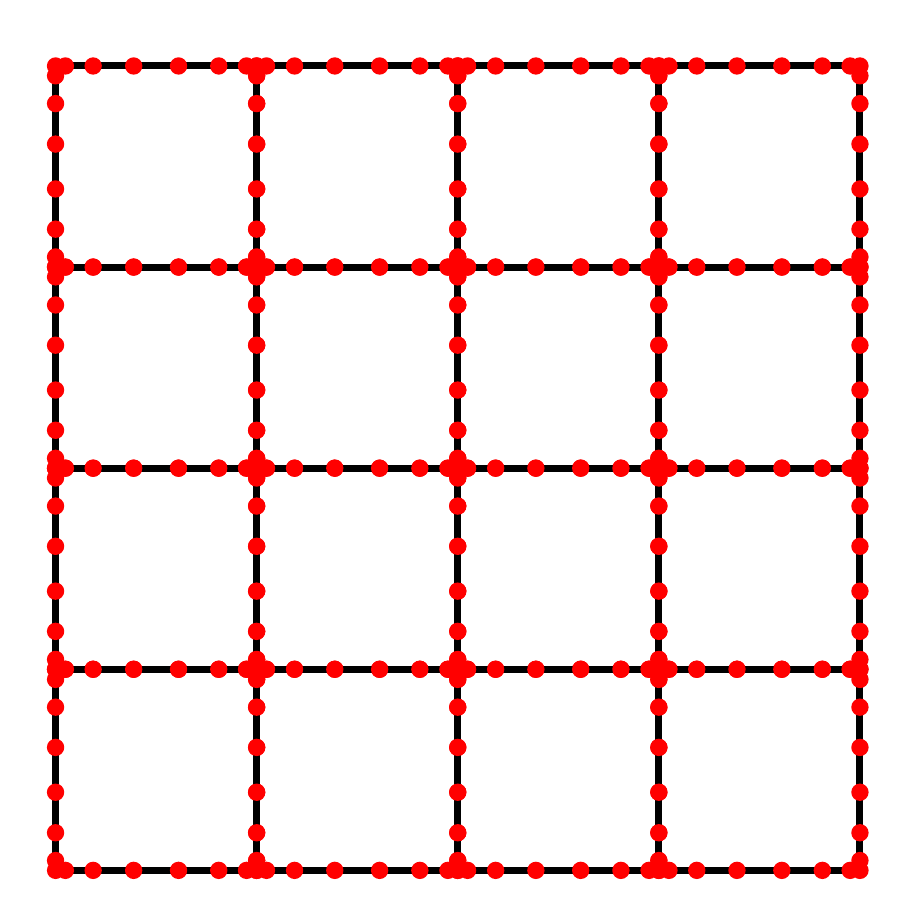}}
\end{picture}
\end{center}

\vspace{2mm}

\noindent
\textbf{Step 2:} Switch tabulation points on the boundary from Chebyshev nodes to
Legendre nodes. The main purpose is to remove the corner nodes.

\begin{center}
\setlength{\unitlength}{1mm}
\begin{picture}(130,57)
\put(00,00){\includegraphics[height=55mm]{fig_HPS_gauss_fig3.pdf}}
\put(63,27){$\Rightarrow$}
\put(60,32){Step 1}
\put(73,00){\includegraphics[height=55mm]{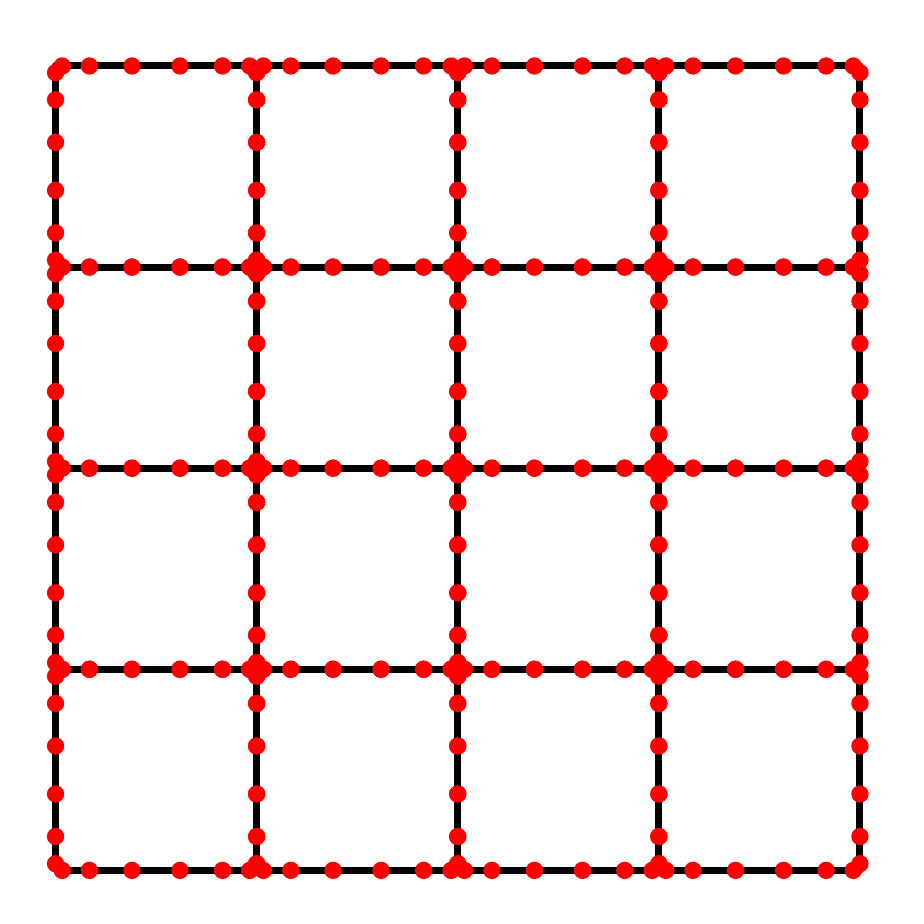}}
\end{picture}
\end{center}

\vspace{2mm}

\noindent
\textbf{Step 3:} Merge the small boxes by pairs as described in Section \ref{sec5:merge}.
The equilibrium equation for each rectangle is formed using the DtN operators of the two
small squares it is made up of. The result is to eliminate the interior nodes (marked in
blue) of the newly formed larger boxes. Construct the solution operator $\UU$ and
the DtN matrices $\VV$ and $\WW$ for the new boxes.

\begin{center}
\setlength{\unitlength}{1mm}
\begin{picture}(130,57)
\put(00,00){\includegraphics[height=55mm]{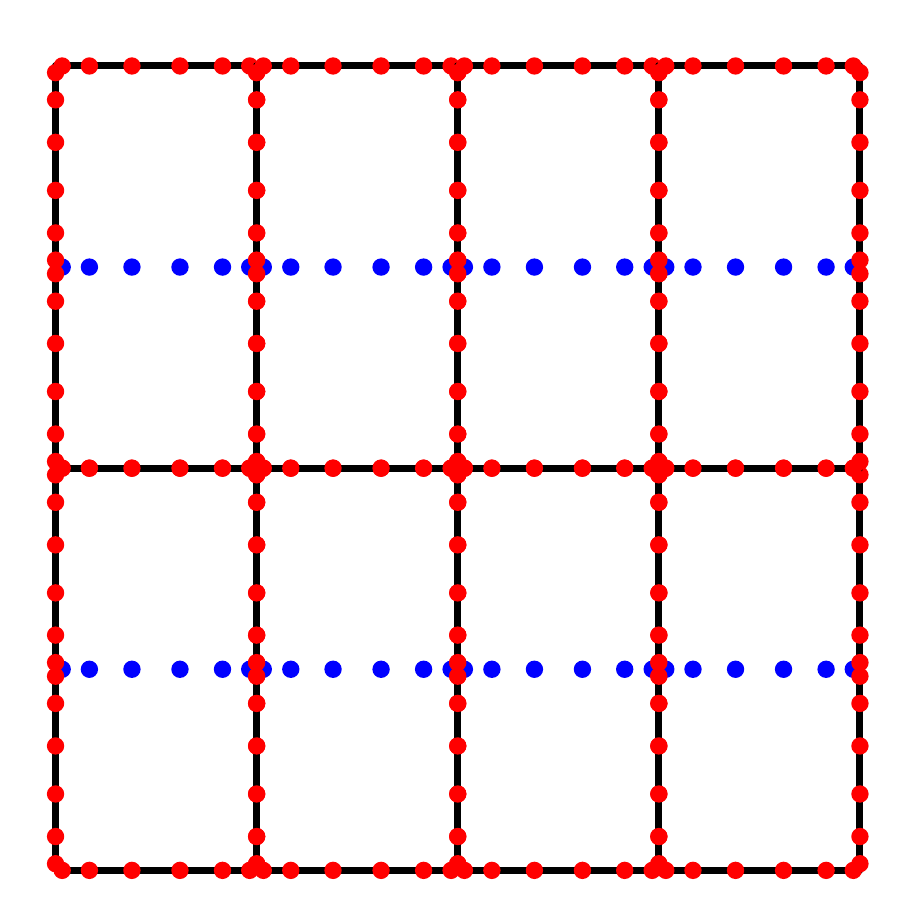}}
\put(63,27){$\Rightarrow$}
\put(60,32){Step 2}
\put(73,00){\includegraphics[height=55mm]{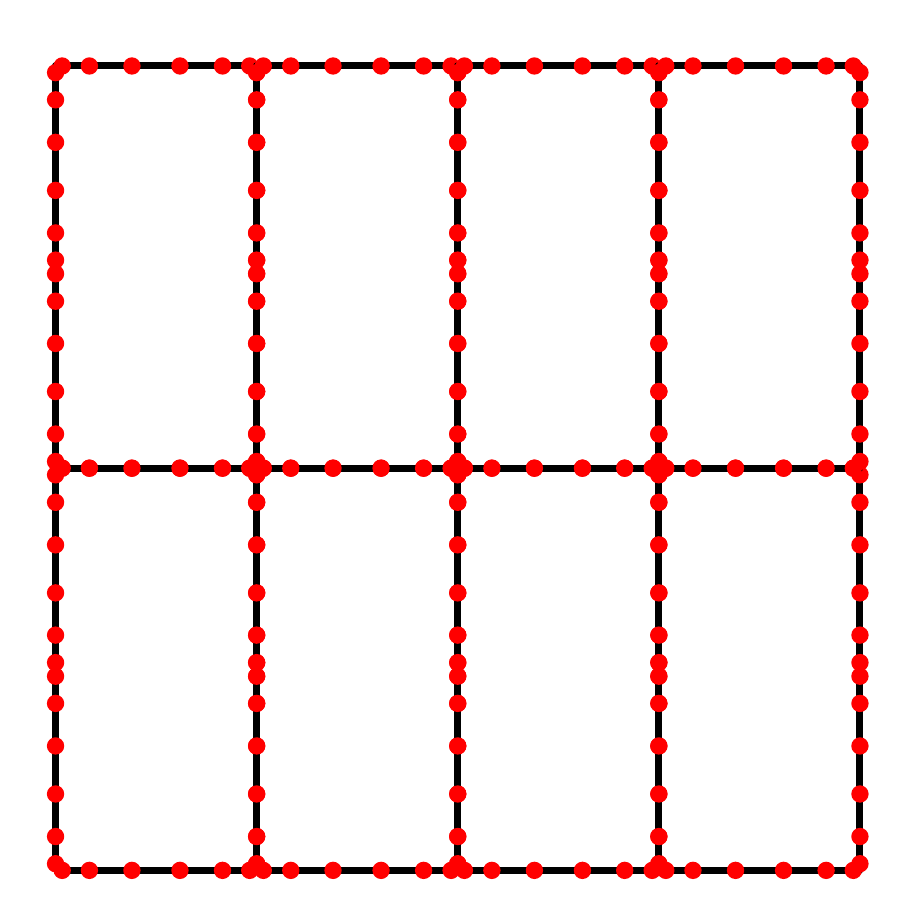}}
\end{picture}
\end{center}

\vspace{2mm}

\noindent
\textbf{Step 4:} Merge the boxes created in Step 3 in pairs, again via
the process described in Section \ref{sec5:merge}.

\begin{center}
\setlength{\unitlength}{1mm}
\begin{picture}(130,57)
\put(00,00){\includegraphics[height=55mm]{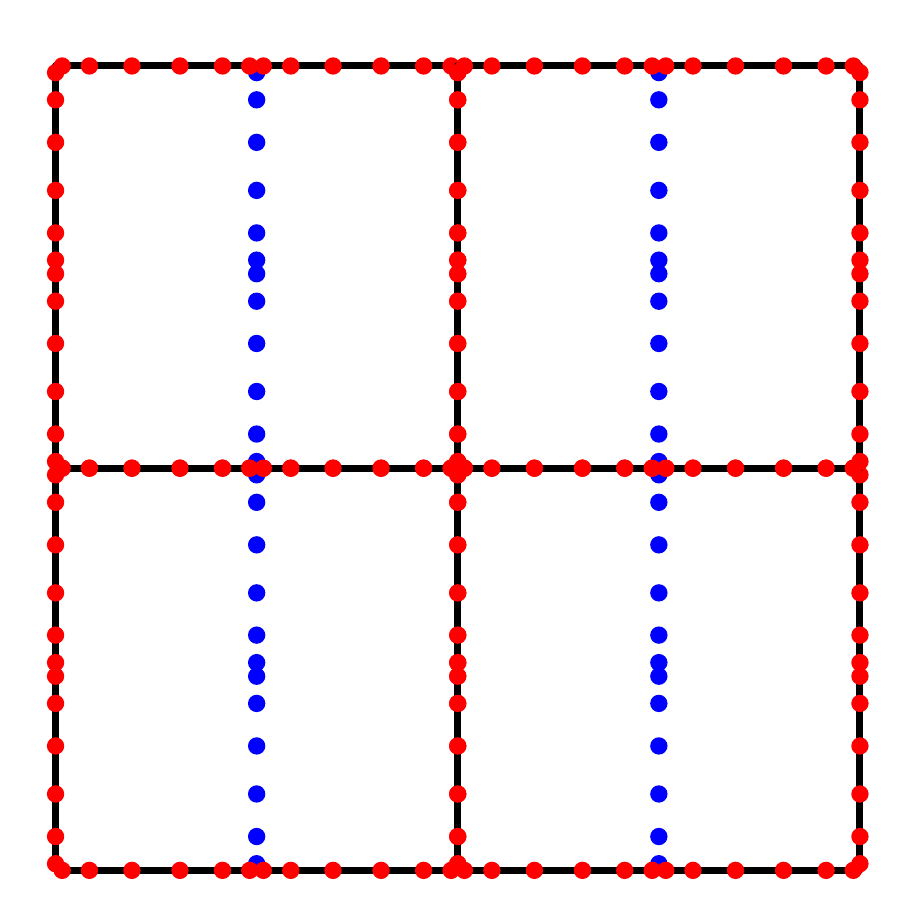}}
\put(63,27){$\Rightarrow$}
\put(60,32){Step 3}
\put(73,00){\includegraphics[height=55mm]{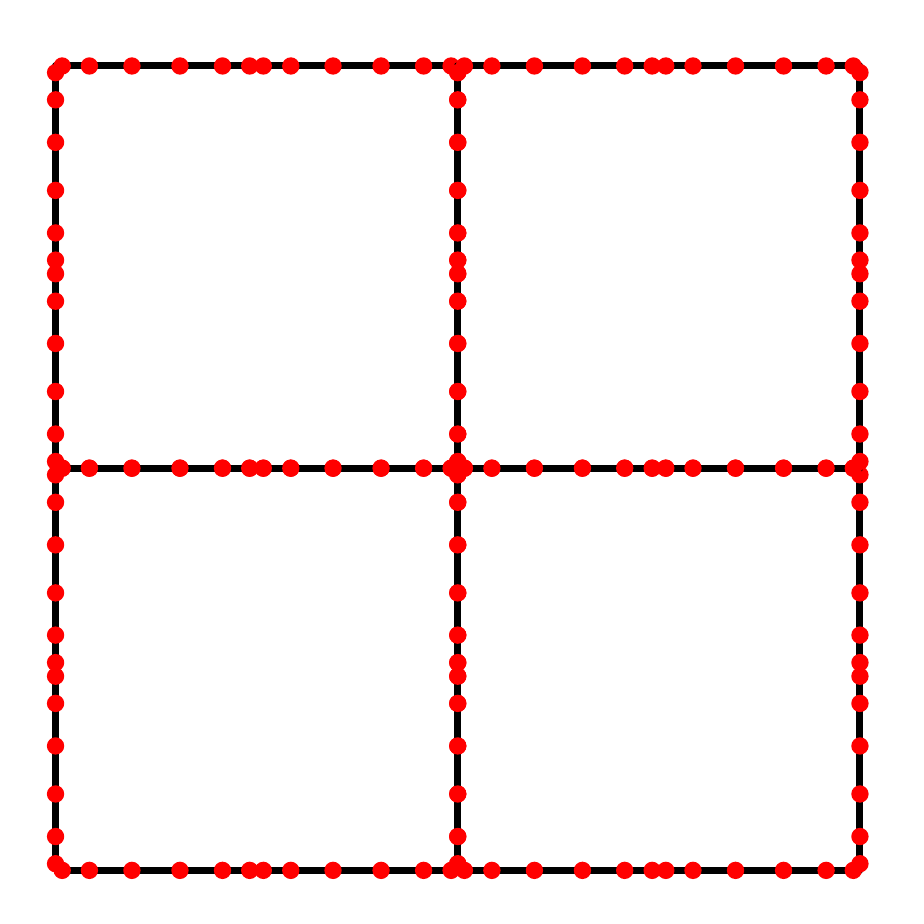}}
\end{picture}
\end{center}

\vspace{2mm}

\noindent
\textbf{Step 5:} Repeat the merge process once more.

\begin{center}
\setlength{\unitlength}{1mm}
\begin{picture}(130,57)
\put(00,00){\includegraphics[height=55mm]{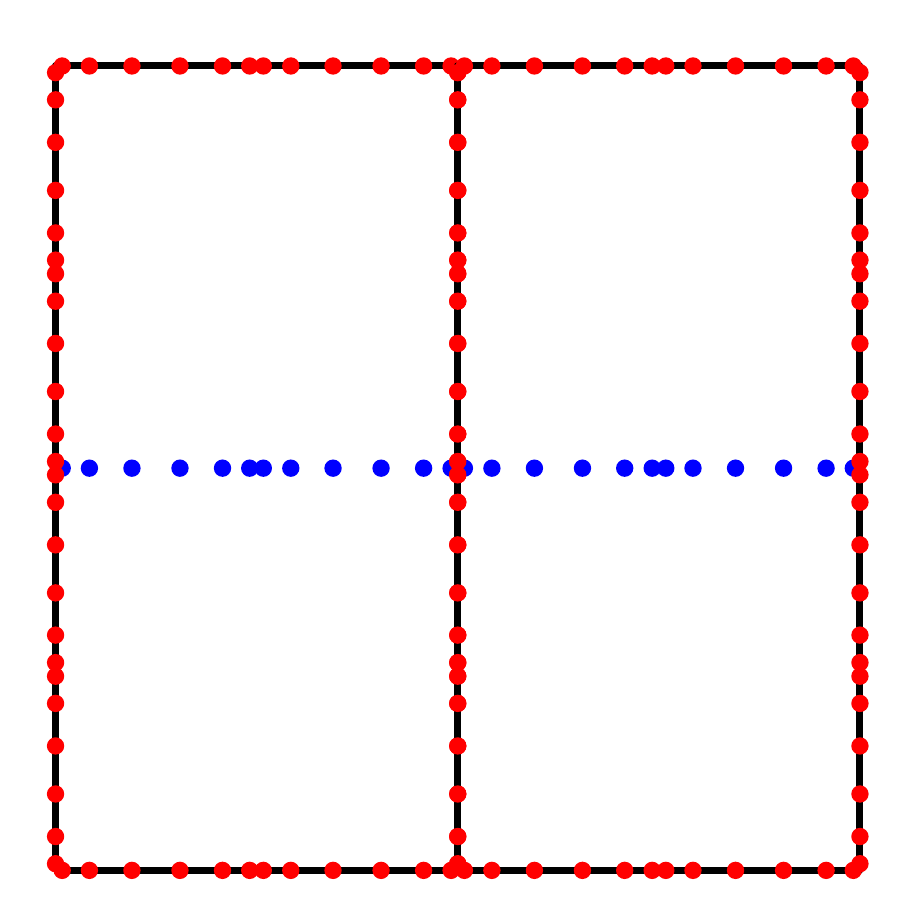}}
\put(63,27){$\Rightarrow$}
\put(60,32){Step 4}
\put(73,00){\includegraphics[height=55mm]{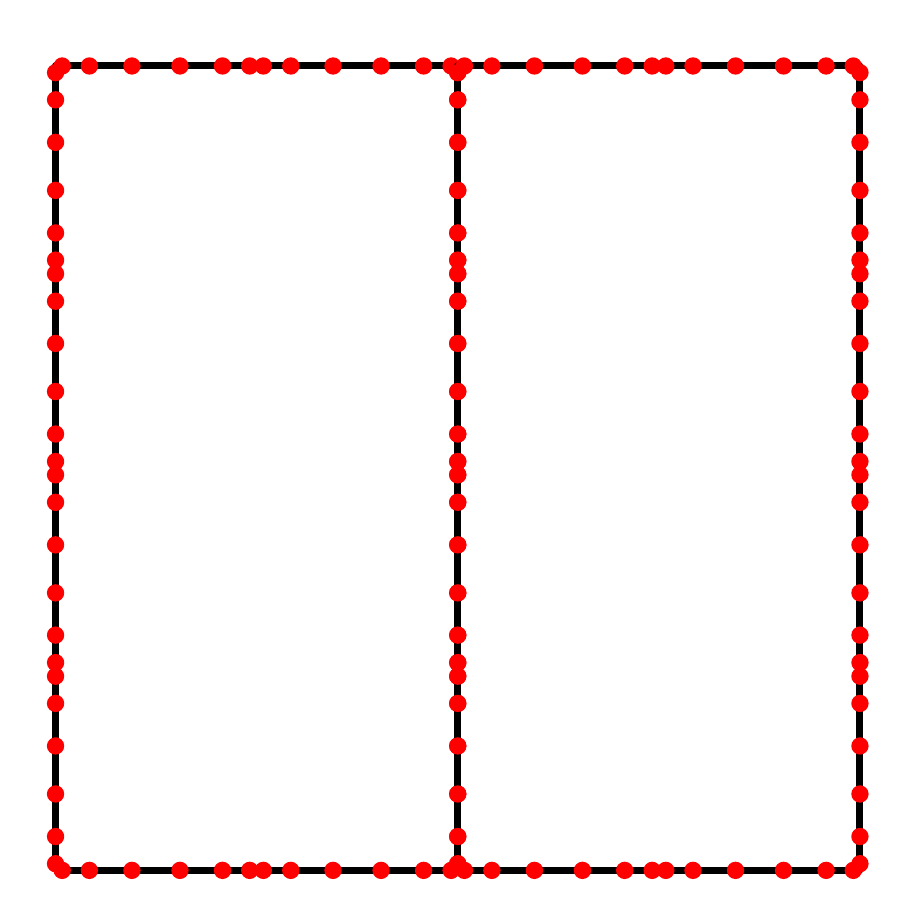}}
\end{picture}
\end{center}

\vspace{2mm}

\noindent
\textbf{Step 6:} Repeat the merge process one final time to obtain the
DtN operator for the boundary of the whole domain.

\begin{center}
\setlength{\unitlength}{1mm}
\begin{picture}(130,57)
\put(00,00){\includegraphics[height=55mm]{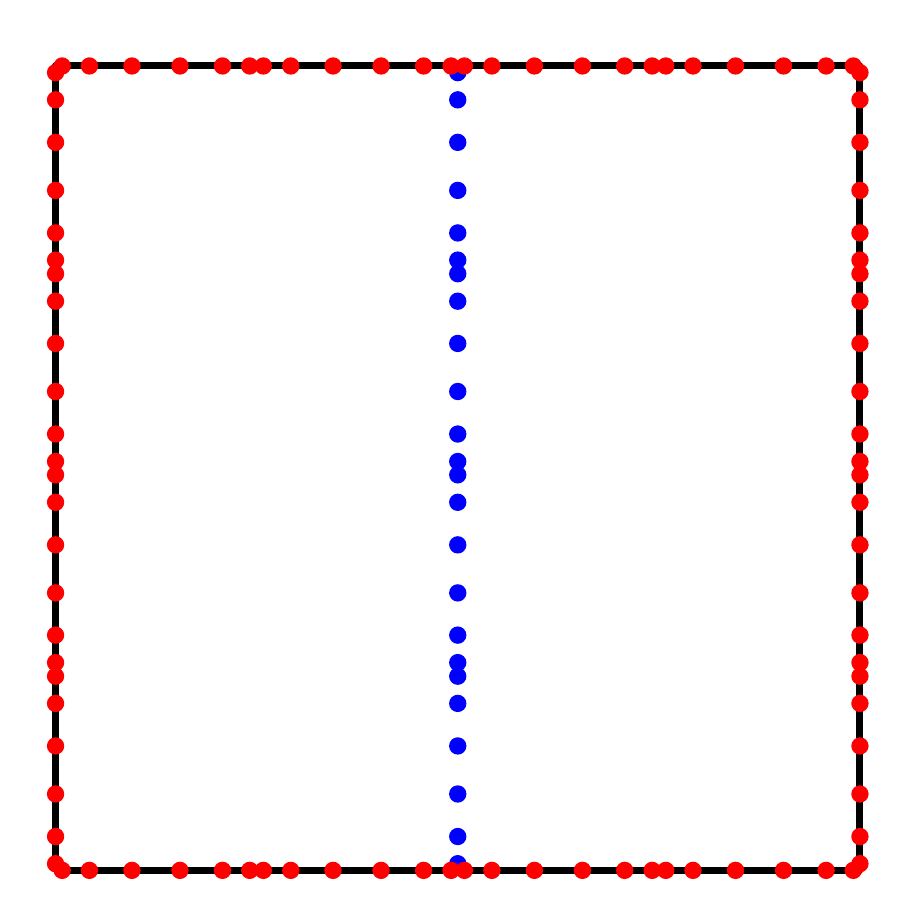}}
\put(63,27){$\Rightarrow$}
\put(60,32){Step 5}
\put(73,00){\includegraphics[height=55mm]{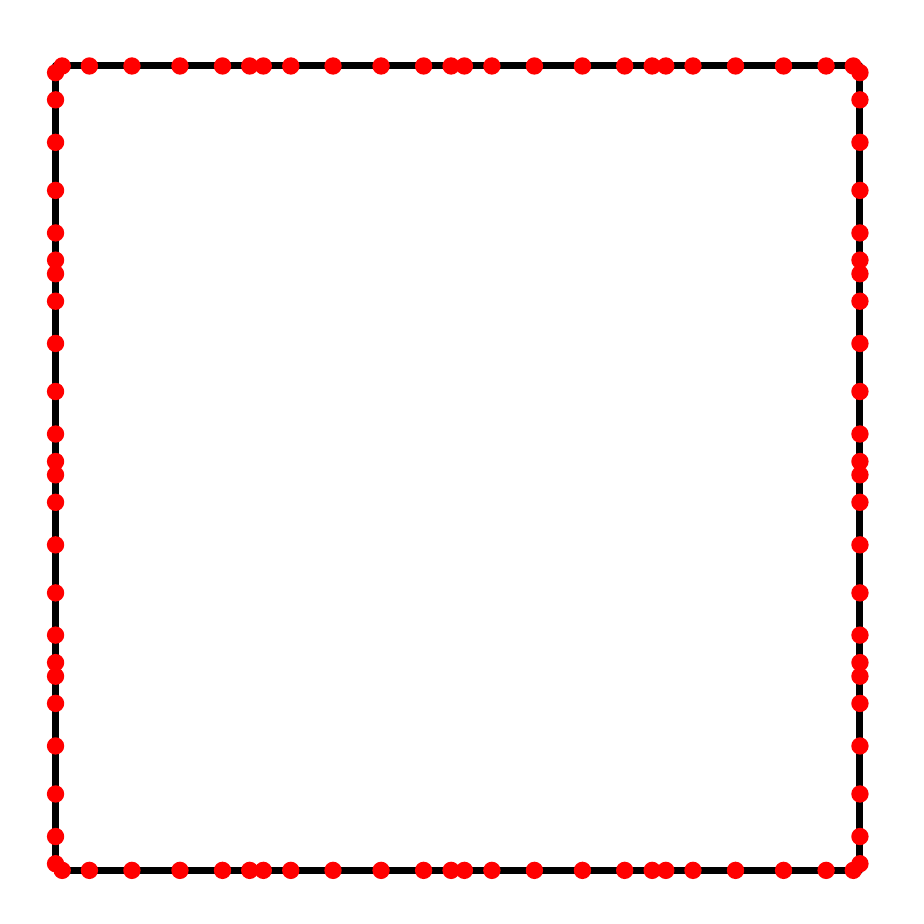}}
\end{picture}
\end{center}

\end{appendix}

\bibliography{main_bib}

\providecommand{\bysame}{\leavevmode\hbox to3em{\hrulefill}\thinspace}
\providecommand{\MR}{\relax\ifhmode\unskip\space\fi MR }
\providecommand{\MRhref}[2]{%
  \href{http://www.ams.org/mathscinet-getitem?mr=#1}{#2}
}
\providecommand{\href}[2]{#2}
\begin{thebibliography}{10}

\bibitem{2008_bebendorf_book}
Mario Bebendorf, \emph{Hierarchical matrices}, Lecture Notes in Computational
  Science and Engineering, vol.~63, Springer-Verlag, Berlin, 2008, A means to
  efficiently solve elliptic boundary value problems. \MR{2451321
  (2009k:15001)}

\bibitem{2010_borm_book}
Steffen B{\"o}rm, \emph{Efficient numerical methods for non-local operators},
  EMS Tracts in Mathematics, vol.~14, European Mathematical Society (EMS),
  Z\"urich, 2010, ${\mathcal{H}}{^{2}}$-matrix compression, algorithms and
  analysis. \MR{2767920}

\bibitem{2006_davis_directsolverbook}
Timothy~A Davis, \emph{Direct methods for sparse linear systems}, vol.~2, Siam,
  2006.

\bibitem{1989_directbook_duff}
I.S. Duff, A.M. Erisman, and J.K. Reid, \emph{Direct methods for sparse
  matrices}, Oxford, 1989.

\bibitem{george_1973}
A.~George, \emph{Nested dissection of a regular finite element mesh}, SIAM J.
  on Numerical Analysis \textbf{10} (1973), 345--363.

\bibitem{2013_martinsson_DtN_linearcomplexity}
A.~Gillman and P.~Martinsson, \emph{A direct solver with $o(n)$ complexity for
  variable coefficient elliptic pdes discretized via a high-order composite
  spectral collocation method}, SIAM Journal on Scientific Computing
  \textbf{36} (2014), no.~4, A2023--A2046, arXiv.org report \#1307.2665.

\bibitem{2013_martinsson_ItI}
Adrianna Gillman, AlexH. Barnett, and Per-Gunnar Martinsson, \emph{A spectrally
  accurate direct solution technique for frequency-domain scattering problems
  with variable media}, BIT Numerical Mathematics \textbf{55} (2015), no.~1,
  141--170 (English).

\bibitem{2002_hackbusch_H2}
W.~Hackbusch, B.~Khoromskij, and S.~Sauter, \emph{On
  $\mathcal{H}^{2}$-matrices}, Lectures on Applied Mathematics, Springer
  Berlin, 2002, pp.~9--29.

\bibitem{hackbusch}
Wolfgang Hackbusch, \emph{A sparse matrix arithmetic based on {H}-matrices;
  {P}art {I}: {I}ntroduction to {H}-matrices}, Computing \textbf{62} (1999),
  89--108.

\bibitem{2014_haut_hyperbolic}
T.S. Haut, T.~Babb, P.G. Martinsson, and B.A. Wingate, \emph{A high-order
  scheme for solving wave propagation problems via the direct construction of
  an approximate time-evolution operator}, arXiv preprint arXiv:1402.5168
  (2014).

\bibitem{1999_hesthaven_spectralcollocation}
J.S. Hesthaven, P.G. Dinesen, and J.P. Lynov, \emph{Spectral collocation
  time-domain modeling of diffractive optical elements}, Journal of
  Computational Physics \textbf{155} (1999), no.~2, 287 -- 306.

\bibitem{2012_martinsson_composite_orig}
P.G. Martinsson, \emph{A composite spectral scheme for variable coefficient
  helmholtz problems}, arXiv preprint arXiv:1206.4136 (2012).

\bibitem{2012_spectralcomposite}
P.G. Martinsson, \emph{A direct solver for variable coefficient elliptic pdes
  discretized via a composite spectral collocation method}, Journal of
  Computational Physics \textbf{242} (2013), no.~0, 460 -- 479.

\bibitem{2003_pfeiffer_spectralmultidomain}
H.P. Pfeiffer, L.E. Kidder, M.A. Scheel, and S.A. Teukolsky, \emph{A
  multidomain spectral method for solving elliptic equations}, Computer physics
  communications \textbf{152} (2003), no.~3, 253--273.

\bibitem{2000_trefethen_spectral_matlab}
L.N. Trefethen, \emph{Spectral methods in matlab}, SIAM, Philadelphia, 2000.

\end{thebibliography}
\bibliographystyle{amsplain}

\end{document}